\documentclass[11pt,b5paper,twoside,headrule]{amsart}
\usepackage{amssymb}
\usepackage{amscd}
\usepackage{amsmath}
\footskip=10mm\parskip 0.2cm \addtocounter{page}{700}
\newtheorem{theorem}{Theorem}[section]
\newtheorem{prop}{Proposition}[section]
\newtheorem{lemma}{Lemma}
\newtheorem{corollary}{Corollary}
\usepackage[]{graphicx}
\usepackage{amsmath,amssymb,amsthm}
\usepackage{verbatim}

\newcommand{\ds}[1]{$\displaystyle{#1}$}

\newtheorem*{theorem*}{Theorem}

\def\text #1{\hbox{\quad#1\quad}}

\newcommand{\Q}{\mathbb Q}
\newcommand{\Z}{\mathbb Z}
\newcommand{\R}{\mathbb R}
\newcommand{\C}{\mathbb C}

\newcommand{\A}{\mathbb A}

\newcommand{\Pp}{\mathbb P}

\makeatletter
\def\l@section{\@tocline{1}{4pt}{1pc}{}{}}
\def\l@subsection{\@tocline{2}{0pt}{2pc}{5pc}{}}
\makeatother

\begin{document}
\font\tenhtxt=eufm10 scaled \magstep0 \font\tenBbb=msbm10 scaled
\magstep0 \font\tenrm=cmr10 scaled \magstep0 \font\tenbf=cmb10
scaled \magstep0


\def\evenhead{{\protect\centerline{\textsl{\large{Dinakar Ramakrishnan and Jonathan Rogawski}}}\hfill}}

\def\oddhead{{\protect\centerline{\textsl{\large{Average Values of Modular $L$-Series}}}\hfill}}

\pagestyle{myheadings} \markboth{\evenhead}{\oddhead}
\thispagestyle{empty}
\noindent{{\small\rm Pure and Applied Mathematics Quarterly\\ Volume 1, Number 4\\
(\textit{Special Issue: In Memory of \\ Armand Borel, Part 3 of 3})\\
701---735, 2005} \vspace*{1.5cm} \normalsize

\begin{center}
{\bf{\Large  Average Values of Modular $L$-series Via the Relative
Trace Formula}}
\end{center}

\begin{center}
\large{Dinakar Ramakrishnan and Jonathan Rogawski}
\end{center}
\footnotetext{Received August 17, 2005. }

\begin{center}\textit {In memory of Armand Borel}
\end{center}

\tableofcontents

\section {\bf Introduction} For any pair of positive
integers $(N,k)$ with $k$ even, let $\mathcal F_N(k)$ denote an
orthogonal basis of holomorphic Hecke-eigencuspforms $\varphi$ of
weight $k$ for the congruence subgroup $\Gamma_0(N)$, acting on the
upper half plane $\mathcal H$. Let $\mathcal F_N(k)^{\rm new}$ be
the subspace of newforms. Denote by $\langle\varphi,\varphi\rangle$
the Petersson norm $\int_{\Gamma_0(N) \backslash \mathcal H}
|\varphi(z)|^2 y^{k-2}dxdy$. Normalize the $L$-function of $\varphi$
so that the functional equation relates $s$ to $1-s$.

For any prime $p$ define probability measures $\mu_+(x)dx$ and
$\mu_-(x)dx$ on the interval $[-2,2] \subset \R$ by
$$
\mu_+(x) \, = \, \frac{(p-1)}{2\pi}
\frac{\sqrt{4-x^2}}{(p^{1/2}+p^{-1/2}-x)^2}
$$
and
$$
\mu_-(x) \, = \, \frac{(p+1)}{2\pi}
\frac{\sqrt{4-x^2}}{\left((p^{1/2}+p^{-1/2})^2-x^2\right)}
$$
Note that $\mu_-$ is the familiar Plancherel measure at $p$, but
(the more interesting) $\mu_+$ is neither that nor the Sato-Tate
measure. In either case, the limit of $\mu_\pm$ exists as $p \to
\infty$ and equals the Sato-Tate measure
$\frac{\sqrt{4-x^2}}{2\pi}dx$.

Now fix a primitive quadratic character $\chi$ of conductor $D<0$,
and restrict attention to the set of $N$ which are prime, do not
divide $D$, and satisfy $\chi(-N) = 1$. From here on assume that $p$
does not divide $DN$, and let $\mu_p(x)$ be $\mu_+(x)$ or $\mu_-(x)$
depending on whether $\chi(p)$ is $1$ or $-1$.

Our main result is the following

\medskip

\noindent{\bf Theorem A} \, \it Let $\chi=\chi_D$ be a quadratic
character of conductor $D<0$, $N$ a prime with $\chi(-N)=1$, $k$
an even integer $>2$, and $p$ an auxiliary prime $\nmid ND$ with
$\chi(p)=\pm 1$. Then for any subinterval $J$ of $[-2,2]$, we have
($\forall \varepsilon>0$)
$$
\sum\limits_{\begin{matrix}\varphi \in \mathcal F_N(k)^{\rm new}
\\ a_p(\varphi) \in J\end{matrix}}
\frac{1}{(\varphi,\varphi)}L(\frac{1}{2}, \varphi \otimes
\chi)L(\frac{1}{2}, \varphi) \, = \, 2\mu_p(J)c_k L(1, \chi) +
O(N^{-k/2+\varepsilon}),
$$
where $L(s, \chi_D)$ is the Dirichlet $L$-function
$\sum\limits_{n\geq 1} \chi_D(n)n^{-s}$, and
\begin{align*}
c_k \, &=\, d({\mathcal D}_k)2^k\pi \frac{k((k/2-1)!)^2}{(k-1)!}\\
&\left(1+\sum_{n=0}^{\frac{k}{2}-2} \left(
\begin{array}{cc}
k\\
2n+1
\end{array}\right)
(-1)^{k/2-n}(k/2+n-1)!(k/2-n-2)!\right),
\end{align*}
with $d({\mathcal D}_k)$ denoting the formal degree of the
holomorphic discrete series representation of PGL$(2, \R)$ of weight
$k$. \rm

\medskip

Here $a_p(\varphi)$ denotes the (normalized) $p$-Hecke eigenvalue of
$\varphi$, which is known by Deligne to lie in $[-2,2]$. If $c_p$ is
the usual $p$-Hecke eigenvalue, which is an algebraic integer, $a_p$
equals $c_p/p^{(k-1)/2}$. Concerning the formal degree, when the
Haar measure $dg = dg_\infty \otimes dg_0$ on GL$(2, \A)$ is chosen
so that $dg_\infty$ corresponds to the Poincar\'e measure $dxdy/y^2$
on the upper half plane, ${\mathcal D}_k = (k-1)/2$.

\medskip

When $J$ is the full interval $[-2,2]$, there is no issue with the
measure, and such a result has then been known for some time (for
$k=2$) by the work of W.~Duke ([Du]); see also [Lu]. Our work began
with our effort to understand these papers. It should be mentioned
that in a recent paper of Iwaniec and Sarnak ([IwS]), they introduce
a novel program to prove a good lower bound for $L(1,\chi)$ from a
simultaneous lower bound, for a family of $\varphi$, of $L(\frac 12
, \varphi)$ and $L(\frac 12 , \varphi\otimes \chi)$. There are other
notable papers dealing with non-vanishing results for modular
$L$-functions and their derivatives; see for example [KMV].

\medskip

Our approach here is representation-theoretic, and makes use of
the {\it relative trace formula} of H.~Jacquet; more precisely, we
study the integral of the GL$(2)$ kernel over the square of the
maximal split torus $T$. In the process we are led to compute and
interpret some period integrals (and certain related local ones)
explicitly, and this is where the technical difficulty resides.
For simplicity of exposition, we are not treating the case $k=2$
here. The reason is that in our method we need to work with a
special, factorizable function $f$ on GL$(2, \A)$ with $f_\infty$
being the matrix coefficient of the discrete series of weight $k$,
which is unfortunately {\it not} integrable for $k=2$. For $k>2$,
this choice of $f_\infty$ in the relative trace formula allows us
to isolate the cusp forms of weight $k$ (and level $N$), and for
$k=2$, we will need to modify $f_\infty$ \`a 1\`a Hecke, which
will force us to consider other cusp forms and Eisenstein series
which have the corresponding $K_\infty$-type. After subtracting
off a suitable singularity on either side of the trace formula, we
will be left with a similar asymptotic formula for $k=2$ as well.

\medskip

Here is a concrete consequence of our theorem:

\medskip

\noindent{\bf Corollary B} \, \it Fix $p, \chi$ as above, with $D$
being the conductor of $\chi$, and an even integer $k > 2$. Choose
any non-empty interval $J \in [-2,2]$. Then there are infinitely
many primes $N$ not dividing $pD$, which are inert in
$K=\Q[\sqrt{D}]$, such that for a holomorphic newform $\varphi$ of
level $N$ and weight $k$,
\begin{enumerate}
\item[(i)] $a_p \in J$; \, \, and
\item[(ii)] $L(1/2,\varphi\otimes \chi)L(1/2,\varphi) \ne 0$.
\end{enumerate}
\rm

\medskip

In particular, there are cusp forms $\varphi$ of weight $k$ and
prime level for which $a_p$ lies arbitrarily close to $2$ or $-2$,
but at the same time having $L(1/2,\varphi\otimes
\chi)L(1/2,\varphi)$ stay non-zero. Without the requirement on the
non-vanishing of the $L$-value, this is a well known result of Serre
([Se]); see also Sarnak ([Sa]) in the context of Maass forms.

\medskip

Finally, we note that, unlike in the Sato-Tate and the Plancherel
cases, there is a {\it bias} built into our measure $\mu_+$. It
favors positive $x$.

\medskip

Much of the work on this paper was done five years ago, but the
explicit nature of the measure was worked out only recently. One can
extend the result without much trouble to square-free level $N$ with
$\chi(-N)=1$, and also to analogous situations over totally real
fields, even with varying weights $\geq 2$ at the different infinite
places, but of the same parity. (If $\chi(-N)=-1$, the relevant
$L$-function vanishes at $s=1/2$ and in that case one should
consider the derivative.) However, the general level $N$ case is
more difficult. One can similarly treat the case of Maass forms, but
then the infinity type $\lambda$ will not be fixed, but will lie in
a short interval. We thank Nathaniel Grossman for helpful comments
on the evaluation of an archimedean integral involving the
hypergeometric function. The first author also thanks Bill Duke,
Jeff Hoffstein, Herv\'e Jacquet, Wenzhi Luo, Philippe Mchel, and
S.~Rallis for interesting conversations.  After this paper was
finalized, we heard from Philippe Michel of a similar work of
E.~Royer (Bull. Soc. Math. France  128  (2000),  no. 2, 219--248)
involving the average of the {\it single} $L$-function $L(1/2,
\varphi)$, without using representation theory. Last but not least,
we acknowledge with thanks the support from the National Science
Foundation.

\vskip 0.2in

\section{\bf The geometric side: Isolation of the main term}

\bigskip

\subsection{Regularization and the main integral}

Let $G$ denote GL$(2)/\Q$ with center $Z$, and set
$\widetilde{G}=G/Z.$  Consider the kernel
\[
K(x,y)=\sum_{\gamma \in \widetilde{G}(\Q)}f(x^{-1}\gamma y)
\]
where $f$ is a factorizable, smooth function on $\widetilde{G}(\A)$
with some good properties (see section 2.3). If $T$ denotes the
maximal split torus of $G$, put $\widetilde T = T/Z$. The relative
trace formula which we want to make use of, and which is due to
Jacquet ([Ja]), involves integration of this kernel over the square
of $\widetilde T(\A)/\widetilde T(\Q)$ against a character.
Unfortunately, this integral does not make sense. In order to take
care of the convergence problem, we regularize the situation.

\medskip

First a few words about the normalization of Haar measures. Let
$\A$ denote the ring of adeles of $\Q$, and let $I_\Q = \A^\ast$
be its multiplicative group of ideles. Recall that $\Q^\ast$ is a
discrete subgroup of $\A^\ast$, with the quotient $A^\ast/\Q^\ast$
being identified with $\R_+^\ast\times \prod_q \Z_q^\ast$, where
$q$ runs over all the primes. Choose the measure $dx/\vert x
\vert$ on $\R^\ast$, with $dx$ denoting the Haar measure on $\R$.
Normalize the measure on each $\Q_q^\ast$ so that $\Z_q^\ast$ gets
volume $1$. Then if $I_\Q^1$ denote the kernel of norm on $I_\Q$,
the compact group $I_\Q^1/\Q^\ast$ gets volume $1$ under the
quotient measure. Finally, we will normalize, for convenience, the
measures on $K_q$  and $Z_q = Z(\Q_q)$ (for every  prime $q$) such
that $K_qZ_q/Z_q$ gets volume $1$.

\medskip

Now let  $\chi $ be a non-trivial, quadratic idele class character
of $\Q$. For $s_1, s_2 \in \C$, set
\[
I_c(f; s_1,s_2) \, = \, \int_{c}\ \int_{c}\ \ K(\left(
\begin{array}{cc}
a & 0 \\
0 & 1
\end{array}
\right) ,\left(
\begin{array}{cc}
b & 0 \\
0 & 1
\end{array}
\right) )\ \chi (a)|a|^{s_{1}}|b|^{s_{2}}\ d^{*}a\ d^{*}b.
\]
where the subscript $c$ indicates that the integrals are taken over
$a,b\in \Q^\ast\backslash \A^\ast$ such that $c^{-1}<|a|, |b|<c$.

\medskip

We are interested in the following limits:

\[
I(f; s_1,s_2)): = \, \lim_{c\rightarrow \infty }\ I_c(f; s_1,s_2).
\]
and
\[
I(f): = \, \lim\limits_{s_1, s_2 \to 0} I(f; s_1, s_2)
\]

\medskip

Let $T$ be the diagonal subgroup. Set
\[
\widetilde{T}=\{t_{a}\}\ \ \ \ \ \ \ \text{where }\ \ \ \ \ \ \ t_{a}=\left(
\begin{array}{cc}
a & 0 \\
0 & 1
\end{array}
\right) .
\]
Clearly, $\widetilde T$ identifies with $T/Z \, \subset G/Z$.

\medskip

For $\delta \in \widetilde{G}$, define the subgroup
\[
C_{\delta }=\{(t,t^{\prime })\in \widetilde{T}\times \widetilde{T}%
:t^{-1}\delta t^{\prime }=z\delta \text{ for some }z\in Z\}.
\]
Then
\[
I_c(f; s_1,s_2) \, = \, \sum_{\{\delta \}}\ \ I_{c}(\delta, f;
s_1,s_2)
\]
where $\{\delta \}$ is a set of representatives for the double cosets $%
\widetilde{T}\backslash \widetilde{G}/\widetilde{T}$ and
\[
I_{c}(\delta, f;s_1,s_2) =\int_{C_{\delta }(\Q)\backslash \left( \widetilde{T}(\mathbf{%
A)\times }\widetilde{T}(\A)\right) }^{c}\ f(t_{a}^{-1}\delta t_{b})\
\chi(a)|a|^{s_{1}}|b|^{s_{2}}\ d^{*}a\ d^{*}b
\]
where the superscript $c$ indicates that the integral is taken over $%
c^{-1}<|a|,\ \ |b|<c$.

\bigskip

\subsection{Coset representatives and centralizers}

\medskip

Define matrices
\[
\xi (x)=\left(
\begin{array}{cc}
1 & x \\
1 & 1
\end{array}
\right)
\]
and
\[
n^{+}=\left(
\begin{array}{cc}
1 & 1 \\
0 & 1
\end{array}
\right) ,\ \ n^{-}=\left(
\begin{array}{cc}
1 & 0 \\
1 & 1
\end{array}
\right) ,\ \ \varepsilon =\left(
\begin{array}{cc}
0 & 1 \\
1 & 0
\end{array}
\right)
\]
Further let $e$ denote the identity matrix in $G$.

Observe that $T\backslash G/T=\widetilde{T}\backslash \widetilde{G}/%
\widetilde{T}.$

\medskip

\begin{lemma}
The set of matrices
\[
\left\{ \xi (x):x\ne 0,1\right\} \cup \left\{ e,\ \ \varepsilon ,\ \ n^{+},\
\ \varepsilon n^{+},\ \ n^{-},\ \ \varepsilon n^{-}\right\}
\]
is a set of representatives for the double cosets $T\backslash G/T.$
\end{lemma}

\medskip

The elements $\xi (x)$ and the orbits they represent will be called \emph{%
regular. }The six remaining representatives and their orbits will be called
\emph{singular}.

We have
\[
\left(
\begin{array}{cc}
a^{-1} &  \\
& 1
\end{array}
\right) \left(
\begin{array}{cc}
x & y \\
z & w
\end{array}
\right) \left(
\begin{array}{cc}
b &  \\
& 1
\end{array}
\right) =\left(
\begin{array}{cc}
ba^{-1}x & a^{-1}y \\
bz & w
\end{array}
\right) .
\]

\medskip

\begin{lemma}
$C_{\delta }=\{e\}$ if $\delta $ is regular or if $\delta \in \left\{ \
n^{+},\ \ \varepsilon n^{+},\ \ n^{-},\ \ \varepsilon n^{-}\right\} .$ On
the other hand,
\begin{eqnarray*}
C_{e} &=&\{(t_{a},t_{a}):a\in I_{\Q}\} \\
C_{\varepsilon } &=&\{(t_{a},t_{a}^{-1}):a\in I_{\Q}\}.
\end{eqnarray*}
\end{lemma}

\medskip

The proof is straight-forward and will be left to the reader.

\medskip

\begin{lemma}
As $\chi$ is non-trivial, we have ($\forall c>0$)
$$
I_{c}(e,f) \, = \,  I_c(\varepsilon, f) \, = \, 0.
$$
\end{lemma}

\medskip

\proof
To treat $\delta =e,$ write a typical element of $\widetilde{T}\mathbf{%
\times }\widetilde{T}$ as $(t_{x}t_{a},t_{a}).$ Then we may write
$I_{c}(e,f)$ as a double integral
\[
\int_{\Q^\ast\backslash\A^\ast, \,\ c^{-1}<|a|<c}\ \ \left( \int_
{x\in\A^\ast, \, c^{-1}<xa<c} f(t_{x}^{-1})\chi(x)\ |x|^{s_{1}}\
d^{*}x\right) \chi(a)\ |a|^{s_{1}+s_{2}}d^{*}a.
\]
Since $\chi$ is non-trivial and of finite order, its restriction to
$I_{\Q}^{1}$ is also non-trivial. Consequently its integral over
$a\in \Q^{*}\backslash I_{\Q}^{1}$ vanishes, showing that
$I_{c}(e,f) = 0.$

\medskip

To treat $\delta =\varepsilon ,$ write a typical element of $\widetilde{T}%
\mathbf{\times }\widetilde{T}$ as $(t_{x}t_{a},t_{a}^{-1}).$ Then we
may write $I_{c}(\varepsilon,f)$ as a double integral
\[
\int_{\Q^\ast\backslash\A^\ast, \,\ c<|a|<c^{-1}}\ \ \left( \int_
{x\in\A^\ast, \, c^{-1}<xa<c} f(t_{x}^{-1})\chi(x)\ |x|^{s_{1}}\
d^{*}x\right) \chi(a)\ |a|^{s_{1}-s_{2}}d^{*}a.
\]
Again, the integral of $\chi$ over $a\in \Q^{*}\backslash
I_{\Q}^{1}$ vanishes, implying that $I_{c}(\varepsilon, f))=0.$

\endproof

\bigskip

\subsection{The test function}

\medskip

Let $S^{\prime }$ be the set of finite primes $q$ at which $\chi$
is ramified. Fix two distinct primes $p ,N\notin S^{\prime }$ and
set
$$
S \, = \, S^{\prime }\cup \{p,N, \infty\}.
$$
We take our function to be of the form
\[
f=f_{\infty }\times f_{p }\times f_{N}\times f_{S^{\prime }}\times
f^{S}
\]
where $f_{S^{\prime }}=\prod_{v\in S^{\prime }}f_{v}$ and $%
f^{S}=\prod_{v\notin S}f_{v}$.

\medskip

Let ${\mathcal D} _{k}$ be the holomorphic discrete series
representation of PGL$(2, \R)$ of lowest weight $k$. We choose the
archimedean function $f_{\infty }$ to be the formal degree
$d({\mathcal D}_k)$ times the complex conjugate of the matrix
coefficient
$$
\left\langle {\mathcal D} _{k}(g)v,v\right\rangle,
$$
where $v$ is a lowest vector of weight $k$, which is unique up to
a scalar. We will take $v$ to be a unit vector. Explicitly we take
\[
f_{\infty }(g)=\ d({\mathcal D}_k)\frac{\left( 2\sqrt{\det g}\right)
^k} {\left(a+d+i(b-c)\right)^{k}} \ \ \ \ \ \ \text{if }g=\left(
\begin{array}{cc}
a & b \\
c & d
\end{array}
\right) \text{ and }\det g>0\text{ }\
\]
and $f_{\infty }(g)=0$ if $\det g<0.$ This function takes the value
$d({\mathcal D}_k)$ at the identity matrix and has the right
$K_\infty$-type, as needed.

\medskip

We take the function $f_{p }$ to be an arbitrary element in the
bi-$Z_{p}K_{p}$-invariant Hecke algebra of $G(\Q_{p}).$ The
function $f_{N}$ is the characteristic function of the subgroup
$Z_{N}K_{0}(N)$ divided by the measure $V_N$ of
$Z_{N}K_{0}(N)/Z_{N},$ where
\[
K_{0}(N)=\left\{ \left(
\begin{array}{cc}
a & b \\
c & d
\end{array}
\right) \in GL_{2}(\mathbf{Z}_{N}):c\equiv 0 \, (\bmod \, N)\right\}
.
\]
We will choose the function $f^{S}$ to be the characteristic
function of $\prod_{v\notin S}Z_{v}GL_{2}(\mathbf{Z}_{v}).$

Finally let $v = q$ be in $S^{\prime}$. If $q^m$ is the conductor of
$\chi_q$, denote by $X$ the set of representatives in $q^{-m}\Z_q
\subset \Q_q$ for $q^{-m}\Z_q/\Z_q$, which is a finite group
isomorphic to $\Z/q^m$. We may view $\chi_v$ as a character of $X$.
Put
$$
f_v \, = \, g(\chi_v)^{-1}\sum_{z \in X} \chi_{1,v}(z) f_{z,v},
$$
with $f_{z,v}$ being the characteristic function of $\left(
\begin{array}{cc}
1 & z \\
0 & 1
\end{array}
\right) K_vZ_v$, and $g(\chi_v)$ the {\it Gauss sum}
$$
g(\chi_v) \, = \, \int_{\Z_q^\ast} \chi_v(x)\psi_v(q^{-m}x) d^\ast
x.
$$
Here $\psi$ denotes the additive character of $\Q_q$ defined as the
composite
$$
\Q_q \, \rightarrow \, \Q_q/\Z_q \, \rightarrow \, \Q/\Z \,
\rightarrow \, S^1,
$$
with the last arrow on the right being $x \to e^{2\pi ix}$. It is
well known that $g(\chi_v)$ has absolute value $q^{m/2}$. Similarly,
the global Gauss sum $g(\chi)$, which is a product of local ones,
has absolute value $\vert D\vert^{1/2}$ since $\vert D\vert$ is the
conductor of $\chi$. In fact we have, due to the oddness of $\chi$,
$$
g(\chi) \, = \, i\vert D\vert^{1/2},
$$
which is the unique square-root of $D$ in $\C$ with positive
imaginary part.

\bigskip

The integral $I(\delta,f )$ factors as a product:

\[
I(\delta,f )\ =\ I_{\infty }(\delta ,f)\ I_{p}(\delta, f )\
I_{N}(\delta ,f)\ I_{S^{\prime }}(\delta ,f)\ I^{S}(\delta ,f),
\]

\medskip

\noindent where $I^S(\delta, f)$ (resp. $I_{S^\prime}(\delta, f)$)
is the product of $I_v(\delta, f)$ over places $v \notin S$ (resp.
$v \in S^\prime$).

\bigskip

\subsection{The dominant terms}

\medskip

It will turn out that since $\chi$ is non-trivial, the terms
corresponding to $\delta = n^\pm$ will be the dominant ones when the
level becomes large. We will first evaluate $I(n^+,f)$, which is a
limit of $I(n^+,f;s_1,s_2)$ as $s_1, s_2 \to 0$.

\medskip

\begin{theorem}
We have
$$
I(n^+,f) \, = \, \frac{\pi k}{2^kd({\mathcal D}_k)\pi
V_N}\left(\sum\limits_{n=0}^{k/2-2}
\begin{pmatrix}k+2n-1\\2n+1\end{pmatrix}\right)F_p(\chi)L(1,\chi),
$$
where
$$
F_p(\chi): = \, \lim\limits_{s_1,s_2\to 0}\frac{I_p(n^+,f;
s_1,s_2)}{L(s,\chi_p)}
$$
and
$$
V_N: = \, {\rm vol}\left({Z_N K_0(N)/Z_N}\right)
$$
\end{theorem}

\medskip

Note that $L(s, \chi_p)$ has a simple pole at $s=0$ iff $\chi(p)=1$.

\medskip

For any place $v$, the local integral at $v$ can be written as:

\begin{eqnarray*}
I_{v}(n^{+};s_1,s_2) &=&\iint_{\Q_{v}^{*}\times \Q_{v}^{*}}\
f_v(\left(
\begin{array}{cc}
ab & a \\
0 & 1
\end{array}
\right) )\ \chi _{v}(a)^{-1}|a|^{-s_{1}}|b|^{s_{2}}\ d^{*}a\
d^{*}b \\
&=&\iint_{\Q_{v}^{*}\times \Q_{v}^{*}}\ f_v(\left(
\begin{array}{cc}
b & a \\
0 & 1
\end{array}
\right) )\ \chi _{v}(a)|a|^{-s_{1}-s_{2}}|b|^{s_{2}}\ d^{*}a\
d^{*}b.
\end{eqnarray*}

\medskip

Theorem 2.1 is a consequence of the following, thanks to the
functional equation of $L(s,\chi)$, which implies
$L(0,\chi)=\frac{\vert D\vert^{1/2}}{\pi}L(1,\chi)$, and the fact
that $\vert D\vert$ is the global conductor of $\chi$:

\medskip

\begin{prop}
\begin{description}
\item [(a)] Let $v$ be a place $\notin\{p,N,\infty\}$. If $\chi_v$ is unramified
(resp. ramified of conductor $q_v^m$), we have: \,
$$
I_v(n^+;s_1,s_2) \, = \, L_v(-s_1-s_2, \chi_{v}),
$$
(resp. $I_v(n^+;s_1,s_2)=g(\chi_v)^{-1}L_v(-s_1-s_2, \chi_{v})$),
implying that
$$
I(n^+,f) \, = \,
g(\chi)^{-1}F_\infty(\chi)F_N(\chi)F_p(\chi)L(0,\chi),
$$
where for $v=\infty, p, N$,
$$
F_v(\chi) \, = \, \lim\limits_{s_1,s_2\to 0}\frac{I_v(n^+,f;
s_1,s_2)}{L(-s_1-s_2,\chi_v)}
$$
\item[(b)] At $v = N$, we have
$$
I_v(n^+;s_1,s_2) \, = \, \frac{1}{V_N} \, L_N(-s_1-s_2, \chi_{N});
$$
\item [(c)] At the archimedean place, since
$\chi_\infty = {\rm sgn}$, the local factor
$L(s,\chi_\infty)=\pi^{-(s+1)/2}\Gamma(\frac{s+1}{2})$ is regular at
$s=0$ with value $1$, and we have
$$
F_\infty(n^+) \, = \, I_\infty(n^+) \, = \, -\frac{2^ki\pi
k((k/2-1)!)^2h(k)}{(k-1)!},
$$
where
$$
h(k): = \, 1+\sum_{n=0}^{\frac{k}{2}-2} \left(
\begin{array}{cc}
k\\
2n+1
\end{array}\right)
(-1)^{k/2-n}(k/2+n-1)!(k/2-n-2)!
$$
\end{description}
\end{prop}

\bigskip

\noindent{\it Remark}: Note that $F_\infty(n^+)$ is purely
imaginary, and so is the Gauss sum $g(\chi)$, because $\chi$ is odd.
It follows that the global geometric term $I(n^+)$ is real. It is
not hard to see that it is in fact positive.

\proof

(a) \, First let $v=q$ be outside $S = S' \cup \{\ell, N, \infty\}$,
corresponding to a prime $q \ne \ell, N$. By our choice of $f$,
$f_v$ is the characteristic function of $Z_vK_v$. Consequently,
\[
(a,b)\rightarrow f_{v}\left(
\begin{array}{cc}
b & a \\
0 & 1
\end{array}
\right)
\]
is simply the characteristic function of $\mathbf{Z}_{v}\times \
\mathbf{Z}_{v}^{*} $. Also, $\chi$ is unramified at $v$. So we get
$$
I_v(n^+) \, = \, \sum_{n \geq 0} \chi(q^{-n})q^{n(s_1+s_2)}.
$$
This is a geometric series summing to
$(1-\chi(q)^{-1}q^{s_1+s_2})^{-1}$, whence the assertion of part (a)
in this case.

\medskip

Next suppose $v = q$ is in $S'$, with $q^m$ being the conductor of
$\chi_v$. By construction,
$$
I_v(n^+) \, = \, g(\chi_v)^{-1}\sum_{z \in X} I_{z,v}(n^+),
$$
where
$$
I_{z,v}(n^+) \, = \, \chi_{v}(z) \iint_{\Q_{v}^{*}\times
\Q_{v}^{*}}\ f_{z,v}\left(
\begin{array}{cc}
b & a \\
0 & 1
\end{array}
\right) \ \chi _{v}(a^{-1})|a|^{-s_{1}-s_{2}}|b|^{s_{2}}\ d^{*}a\
d^{*}b.
$$
Only those $z$ with invertible images in $q^{-m}\Z_q/\Z_q$ give a
non-zero contribution. So we will restrict our attention to these.

Recall that $f_v$ has support on $\left(
\begin{array}{cc}
1 & z \\
0 & 1
\end{array}
\right)K_vZ_v$, and note that
$$
\lambda \left(
\begin{array}{cc}
1 & -z \\
0 & 1
\end{array}
\right)\left(
\begin{array}{cc}
b & a \\
0 & 1
\end{array}
\right) \, = \, \left(
\begin{array}{cc}
\lambda b & \lambda(a-z) \\
0 & \lambda
\end{array}
\right).
$$
Then it follows that for $f_{z,v}\left(
\begin{array}{cc}
b & a \\
0 & 1
\end{array}
\right)$ to not vanish, we need  $b \in \Z_q^\ast$ and $a -z \in
\Z_q$. Thus $a$ lies in $q^{-m}\Z_q$ and has the same image as $z$
in $q^{-m}\Z_q/\Z_q$. Consequently, since $\chi_{v}$ has conductor
$q^m$ and thus the pullback of a character of
$(q^{-m}\Z_q/\Z_q)^\ast \simeq (\Z/q^m)^\ast$, we must have
$\chi_{v}(z)\chi_{v}(a^{-1}) = 1$. We get
$$
I_{z,v}(n^+) \, = \, {\rm vol}(1 + q^m \Z_q),
$$
when $z$ has invertible image in $\Z/q^m$. The assertion follows
once we note:
\begin{enumerate}
\item[(i)] \, $\varphi(q^m)$vol$(1+q^m \Z_q)$ equals
vol$(\Z_q^\ast)=1$, \, \, and \item[(ii)] \, $L(s, \chi_{v}) = 1$ as
$\chi_{v}$ is ramified.
\end{enumerate}

\medskip

(b) \, Let $v = N$. Recall that $f_N$ is the characteristic function
of $Z_N K_0(N)$ divided by $V_N$. Like in (a), the function
\[
(a,b)\rightarrow f_{N}\left(
\begin{array}{cc}
b & a \\
0 & 1
\end{array}
\right)
\]
is the characteristic function of $\mathbf{Z}_{N}\times \
\mathbf{Z}_{N}^{*}$, but divided by $V_N$. The assertion follows.
(When we consider $n^-$, the situation will be slightly different.)

\medskip

(c) \, Let $v = \infty$. Recall that $\chi_{\infty}(-1) = -1$. By
the definition of $f_\infty$, its value on $\left(
\begin{array}{cc}
b & a \\
0 & 1
\end{array}
\right)$ is $d({\mathcal D}_k)\frac{b^{k/2}}{(b+1-ia)^k}$ (resp.
$0$) when $b$ is positive (resp. negative). Noting that $d^\ast x =
dx/|x|$ on $\R^\ast$. We get
\begin{eqnarray*}
d({\mathcal D}_k)^{-1}I_\infty(n^+;s_1,s_2)  &=&
2^k\int_{-\infty}^{\infty} \int_0^\infty
\frac{b^{k/2+s_2-1}|a|^{-s_1-s_2-1}{\rm sgn}(a)}{(b+1-ia)^k} da db \\
&=&2^k\int_{-\infty}^{\infty} \int_0^\infty \frac{b^{k/2+s_2-1}
|a|^{-s_1-s_2-1}{\rm sgn}(a)((b+1)+ia)^k}{(a^2+(b+1)^2)^k} da db,
\end{eqnarray*}
where $b$ runs over $(0, \infty)$ and $a$ runs over $(-\infty,
\infty)$. Appealing to the binomial expansion
$$
((b+1)+ia)^k \, = \, \sum_{j=0}^k \left(
\begin{array}{cc}
k \\
j
\end{array}
\right) i^{k-j}a^{k-j}(b+1)^j,
$$
we may write
$$
d({\mathcal D}_k)^{-1}I_\infty(n^+;s_1,s_2) \, = \, 2^k\sum_{j=0}^k
\left(
\begin{array}{cc}
k \\
j
\end{array}
\right) i^{k-j} I_j,
$$
where
$$
I_j \, = \, \int_{-\infty}^{\infty} \int_0^\infty
\frac{b^{k/2+s_2-1}|a|^{k-j-s_1-s_2-1}(a/|a|)^{1+k-j}(b+1)^j}{(a^2+(b+1)^2)^k}
da db.
$$

Note that that if $j$, and hence $k-j$, is even, then the
integrand is odd in the $a$-variable and so $I_j$ vanishes. So we
may, and we will, assume from hereon that $j$ is odd. We have
\begin{eqnarray*}
I_j  &=&  2 \int_0^\infty \int_0^\infty \frac{b^{k/2+s_2-1}a^{k-j-s_1-s_2-1}(b+1)^j}{(a^2+(b+1)^2)^k} da db \\
&=&2\int_0^\infty \int_0^\infty \frac{b^{k/2+s_2-1}a^{k-j-s_1-s_2-1}(b+1)^{j-2k}}{((a/(b+1))^2+1)^k} da db \\
&=&2\int_0^\infty \int_0^\infty \frac{b^{k/2+s_2-1}a^{k-j-s_1-s_2-1}(b+1)^{-k-s_1-s_2}}{(a^2+1)^k} da db.
\end{eqnarray*}
So the integral factors as
$$
I_j \, = \, I_{j,1}I_{j,2},
\leqno(0)
$$
where
$$
I_{j,1} \, = \, 2\int_0^\infty \frac{a^{k-j-s_1-s_2-1}}{(a^2+1)^k} da
$$
and
$$
I_{j,2} \, = \, \int_0^\infty b^{k/2 +s_2 -1} (b+1)^{-k-s_1-s_2} db.
$$
Note that $I_{j,2}$ is in fact independent of $j$.

\medskip

By the substitution $u = a^2$, we get
$$
I_{j,1} \, = \, \int_0^\infty \frac{u^{(k-j-s_1-s_2)/2 - 1}}{(u+1)^k} du.
$$
Recall that the Beta function $B(z,w) =
\frac{\Gamma(z)\Gamma(w)}{\Gamma(z+w)}$ has the following integral
representations (cf. [A-S], p.258):
$$
B(z,w) \, = \, \int_0^\infty t^{z-1}(1-t)^{w-1} dt \, = \, \int_0^\infty \frac{t^{z-1}}{(1+t)^{z+w}} dt.
\leqno(*)
$$
It follows immediately that
$$
I_{j,1} \, = \, B((k-j-s_1-s_2)/2, (k+j+s_1+s_2)/2)
\leqno(1)
$$
Now let $v = b/(b+1)$ in the expression for $I_{j,2}$, so that $b =
v/(1-v)$ and $db = dv/(1-v)^2$. We obtain:
$$
I_{j,2} \, = \, \int_0^1 v^{k/2-1+s_2} (1-v)^{k/2 + s_1 -1} dv.
$$
Applying the first identity of (*) we get
$$
I_{j,2} \, = \, B(k/2 + s_2, k/2 + s_1).
\leqno(2)
$$
Putting $(0), (1)$ and $(2)$ together, and writing everything in
terms of the $\Gamma$-function, we finally obtain
$$
I_j \, = \, \frac{\Gamma((k-j-s_1-s_2)/2)\Gamma((k+j+s_1+s_2)/2)
\Gamma(k/2+s_2)\Gamma(k/2+s_1)}{\Gamma(k)\Gamma(k+s_1+s_2)}.
\leqno(3)
$$
Consequently, we see that $d({\mathcal
D}_k)^{-1}I_\infty(n^+;s_1,s_2)$ equals
$$
2^k\left(\sum_{\begin{array}{cc}j=0\\
j \equiv 1 (2)
\end{array}
}^k \left(
\begin{array}{cc}
k\\
j
\end{array}
\right)
i^{k-j}\Gamma(\frac{k-j-s_1-s_2}{2})\Gamma(\frac{k+j+s_1+s_2}{2})\right)
\frac{\Gamma(\frac{k}{2}+s_1)\Gamma(\frac{k}{2}+s_2)}{\Gamma(k)\Gamma(k+s_1+s_2)}.
$$
Now observe that for any positive integer $r$,
$$
\Gamma\left(\frac{1+2r}{2}\right) \, = \, (r-1)!\Gamma(\frac12).
$$
Since $\Gamma(\frac12)=\sqrt{\pi}$ and $\Gamma(r)=(r-1)!$, putting
$k=2m, j=2n+1$ (so that $i^{k-j}=-i(-1)^{m-n}$), we get (by sending
$s_1, s_2$ to $0$):
\begin{align*}
I_\infty(n^+) \, &= \, -2^ki\pi d({\mathcal
D}_k)\frac{(2m)((m-1)!)^2}{(2m-1)!}\\
 &\left(1+ \sum_{n=0}^{m-2}\left(
\begin{array}{cc}
2m\\
2n+1
\end{array}\right)
(-1)^{m-n}(m+n-1)!(m-n-2)!\right)
\end{align*}
 This proves (c). \qed

\bigskip

Next we consider the case $\delta = n^{-}$. In this case, the local integral at any $v$ equals
\begin{eqnarray*}
I_{v}(n^{-};s_1,s_2) &=&\iint_{F_{v}^{*}\times F_{v}^{*}}\
f_v(\left(
\begin{array}{cc}
ab & 0 \\
b & 1
\end{array}
\right) )\ \chi _{v}(a)^{-1}|a|^{-s_{1}}|b|^{s_{2}}\ d^{*}a\
d^{*}b \\
&=&\iint_{\Q_{v}^{*}\times \Q_{v}^{*}}\ f_v(\left(
\begin{array}{cc}
a & 0 \\
b & 1
\end{array}
\right) )\ \chi _{v}(a^{-1}b)|a|^{-s_{1}}|b|^{s_1+s_{2}}\ d^{*}a\
d^{*}b.
\end{eqnarray*}

\medskip

\begin{prop}
\begin{description}
\item [(a)] If $v$ is a place outside $\{p,N,\infty\}$, we have: \,
$$
I_v(n^-;s_1,s_2) \, = \, I_v(n^+:-s_2,-s_1);
$$
\item[(b)] At $v = N$, we have
$$
I_N(n^-;s_1,s_2) \, = \, \frac{1}{V_N} \chi(N)N^{-s_1-s_2}
L(s_1+s_2, \chi_{N});
$$
\item [(c)] At the archimedean place we have:
$$
I_\infty(n^-) = -I_\infty(n^+)
$$
\end{description}
\end{prop}

\medskip

Since $\chi(N)=-1$, the sign change at infinity is compensated for
by the sign change at $N$, when we go from $n^+$ to $n^-$. We get
the following

\medskip

\begin{corollary} \,
\begin{enumerate}
\item
$$
I(n^-,f)F_p(\chi) \, = \, I(n^+,f)F_p^-(\chi),
$$
where
$$
F_p^-(\chi): = \, \lim\limits_{s_1,s_2\to 0}\frac{I_p(n^-,f;
s_1,s_2)}{L(s,\chi_p)}
$$
\item When $f_p$ is the characteristic function of $Z_pK_p$, we have
$$
F_p(\chi) \, = \, F_p^-(\chi) \, = \, 1.
$$
\end{enumerate}
\end{corollary}

\medskip

{\it Proof of Proposition}. \, Part (a) is immediate. At $N$, the
function $F_N$ when restricted to
\[
(a,b)\rightarrow f_{N}\left(
\begin{array}{cc}
a & 0 \\
b & 1
\end{array}
\right)
\]
is simply the characteristic function of $\Z_N^\ast \times N\Z_N$
divided by $V_N$. Since $\chi$ is unramified at $N$,
$\chi(a^{-1}b)$ is just $\chi(b)$. Now part (b) follows by
switching to the variable $b'$ with $b=Nb'$. At $\infty$, we
obtain
$$
I_\infty(n^-;s_1,s_2) \, = \, d({\mathcal
D}_k)2^k\iint\limits_{\R^\ast\times \R^\ast}
f_{\infty}\begin{pmatrix}a&0\\b&1\end{pmatrix}{\rm
sgn}(a^{-1}b)\vert a\vert^{-s_1}\vert b\vert^{s_1+s_2}d^\ast ad^\ast
b,
$$
and by the definition of $f_\infty$, it is zero if $a<0$ and for
$a>0$ it takes the value $d({\mathcal
D}_k)\frac{a^{k/2}}{((a+1)+ib)^k}$. Comparing with what we began
with at $\infty$ for $n^+$, we see immediately (by reversing the
roles of $a$ and $b$) that
$$
I_\infty(n^-;s_1,s_2) \, = \, -I_\infty(n^+; -s_2,-s_1).
$$
The assertion of part (c) follows by sending $s_1, s_2$ to $0$.

\qed

\bigskip

\subsection{The remaining singular terms}

\medskip

We have already seen that under our hypotheses, the global integrals
$I_c(e,f)$ and $I_c(\varepsilon,f)$ are zero for any $c>0$. So
$$
I(e,f) \, = \, I(\varepsilon,f) \, = \, 0.
$$

\medskip

So it remains only to evaluate $I(\varepsilon n^\pm, f)$.

\bigskip

\begin{prop} \, Let $\delta \in \{\varepsilon n^+, \varepsilon n^-\}$. Then $I(\delta,f) = 0$.
\end{prop}

\medskip

\proof As $I(\delta,f) = \prod_v I_v(\delta,f)$, it suffices to show
that $I_N(\varepsilon n^\pm, f) = 0$.

First consider the $\delta = \varepsilon n^+$ case.

\medskip

We have

\[
I_N(\varepsilon n^{+})=\iint_{\Q_N^\ast \times \Q_N^\ast}\
f_N(\left(
\begin{array}{cc}
0 & a \\
b & 1
\end{array}
\right) )\ \chi(a)^{-1}|a|^{-s_{1}}|b|^{s_{2}}\ d^{*}a\ d^{*}b
\]
By definition, $f_N$ is supported on $Z_N K_0(N)$. So $f_N(\left(
\begin{array}{cc}
0 & a \\
b & 1
\end{array}
\right) )$ is non-zero for some $(a,b) \in (\Q_N^\ast)^2$ iff there
is a $\lambda \in \Q_N^\ast$ such that $\lambda, \lambda a \in
\Z_N$, $\lambda b \in N\Z_N$ and det$(\lambda\left(
\begin{array}{cc}
0 & a \\
b & 1
\end{array}
\right) \, = \, -\lambda^2 ab \in \Z_N^\ast$. These conditions
cannot be satisfied simultaneously. Done.

\medskip

Now let $\delta = \varepsilon n^-$. Then we have
\[
I_N(\varepsilon n^{-})=\int_{(\Q_N^\ast)^2}\ f_N(\left(
\begin{array}{cc}
ab & a \\
b & 0
\end{array}
\right) )\ \chi(a)^{-1}|a|^{-s_{1}}|b|^{s_{2}}\ d^{*}a\ d^{*}b
\]
Again this can be non-zero only if we can find, for each $a, b$, a
$\lambda$ in $\Q_N^\ast$ such that $\lambda ab, \lambda a \in \Z_N$,
$\lambda b \in N\Z_N$ and $\lambda^2 ab \in \Z_N^\ast$, which is
impossible.

\qed

\medskip

To recapitulate, four of the six singular terms contribute zero to
the geometric side of the relative trace formula. The remaining
two, which give the dominant terms, are both non-zero, and in
fact, thanks to the functional equation of $L(s,\chi)$, have the
{\it same} limiting expression, which is a non-zero multiple of
$L(1,\chi)$, as $s_1, s_2$ both go to $0$.

\vskip 0.2in

\subsection{The regular terms}

\medskip

Recall that the regular double cosets are represented by the
matrices $\{\xi(x)\}$, with $x \in \mathbb P^1(\Q) - \{0,1,
\infty\}$. By abuse of notation we will write $I(x)$ for
$I(\xi(x))$. At any place $v$, we have by definition,

\[
I_v(x)=\int_{\Q_v^\ast\times \Q_v^\ast}\ f_v\left(
\begin{array}{cc}
ab & ax \\
b & 1
\end{array}
\right) \ \chi _{v}(a)^{-1}|a|^{-s_{1}}|b|^{s_{2}}\ d^{*}a\ d^{*}b
\]

\medskip

\begin{prop}
Let $x \in \Pp^1(\Q) - \{0,1, \infty\}$.
\begin{description}
\item [(a)] {Let $v = q$ be a prime not in $S$. When $v(1-x) > 0$, $I_v(x)$
vanishes. Suppose $v(1-x) = 0$, resp. $v(1-x) < 0$. Then $v(x) \geq
0$, resp. $v(x) = v(1-x)$, and we have:
$$
I_v(x) \, = \, \sum_{n = 0}^{v(x)} \, {\sum_{m=-v(x)}^{-n}}{}^\prime
\, (\overline \chi(q)q^{s_1})^m q^{-ns_2}),
$$
resp.
$$
I_v(x) \, = \, \sum_{n = v(x)}^{0} \, {\sum_{m=0}^{-n}}{}^\prime \,
(\overline \chi(q)q^{s_1})^m q^{-ns_2},
$$
where the prime on the inside sum implies (in either case) that the
summation is over $m$ of the same parity as $n + v(1-x)$. In
particular,
$$
v(x) = v(1-x) = 0 \, \implies \, I_v(x) = {\rm vol}(\Z_v^\ast)^2.
$$}
\item[(b)] {Let $v = q$ be a prime in $S^\prime$, with $q^c$ the
conductor of $\chi_q$. Then if $v(1-x) > c$, $I_v(x,f))$ vanishes.
When $v(1-x) \leq c$, $I_v(x,f)$ is bounded by
$q^{c/2}(v(x)+c)^2{\rm vol}(\Z_v^\ast)^2$.}
\item[(c)] {Let $v =
N$. Then $I_N(x)$ vanishes when $v_N(1-x) \ne  0$. When $v_N(1-x)
= 0$, we must have $v_N(x) \geq 1$ and
$$
I_N(x) \, = \, \sum_{n = 1}^{v_N(x)}
\sum_{m=-v_N(x)}^{-n}{}^\prime \, (\overline \chi(N)N^{s_1})^m
N^{-ns_2}.
$$}
\item[(d)] {Let $v=p$, and $f_p$ the characteristic function of ${ZK\left(
\begin{array}{cc}
  p^r & 0 \\
  0 & p^{r'} \\
\end{array}
\right)K}$ for some integers $r, r'$ such that $r \geq r'$. Then
$I_p(x,f)$ is zero unless $v(1-x) < r+r'$, in which case it is
bounded by $C(f_p)v(x)^2{\rm vol}(\Z_v^\ast)^2$, for a constant
$C(f_p)$.
\item [(e)] Let $v = \infty$. Assume that $x > 0$ and that $\Re s_1, \Re s_2 \in (-k/2, k/2)$.
Then, as $\chi_\infty(-1)=-1$, we have for $1-x$ is positive (resp.
negative):
$$
I_\infty(x) \, = \, (1-x)^{k/2} \left( I_\infty(-1,1,1)
-(-1)^{k}I_\infty(-1,-1,1)\right)
$$
(resp.
$$
I_\infty(x) \, = \, (x-1)^{k/2} \left( I_\infty(1,-1,-1)
-(-1)^{k}I_\infty(1,1,-1)\right),
$$
where for $\epsilon, \delta, \nu \in \{\pm 1\}$,
$$
I_\infty(\epsilon, \delta, \nu) \, = \, (-\epsilon)^{\rho-\sigma
-1}\delta^{\rho+k-3\sigma-1}\nu^{k-\sigma}i^{\rho-2\sigma}
B(\sigma,k-\sigma)B(\rho, k-\rho)F(k-\sigma,\rho;k;1-\epsilon\nu),
$$
where $F = {}_{2}F_1$ is the hypergeometric function.}
\end{description}
\end{prop}

\medskip

\proof \, (a) \, Since $v = q$ is a prime outside $S$, $f_v$ is by
definition the characteristic function of $Z_vK_v$.

\medskip

\begin{lemma} \, $f_v\left(
\begin{array}{cc}
ab & ax \\
b & 1
\end{array}
\right)$ is non-zero iff the following conditions are
simultaneously satisfied:
\begin{enumerate}
\item[(i)] $v(1-x) \leq 0$;
\item[(ii)] $v(x) \geq v(1-x)$; in particular, $v(x) = v(1-x)$ when $v(1-x) < 0$;
\item[(iii)] $v(1-x)-v(x) \leq v(a) \leq$ min$\{-v(1-x), -v(b)-v(1-x)\}$;
\item[(iv)] max$\{v(1-x), v(a)+v(1-x)\} \leq v(b) \leq v(x)-v(1-x)$;
\item[(v)] $v(a) + v(b) \equiv v(1-x)$ \, $($mod $2)$.
\end{enumerate}
\end{lemma}

\medskip

{\it Proof of Lemma}. \, For $f_v\left(
\begin{array}{cc}
ab & ax \\
b & 1
\end{array}
\right)$ to be non-zero, it is necessary and sufficient that there
exists a $\lambda$ in $\Q_v^\ast$ such that $\lambda ab, \lambda
ax, \lambda b, \lambda$ are in $\Z_v$ {\it and} $\lambda^2
ab(1-x)$ is in $\Z_v^\ast$; in other words, one must have
$$
2v(\lambda) + v(a) + v(b) + v(1-x) = 0;
\leqno 1)
$$
$$
v(\lambda) + v(a) + v(b) \geq 0;
\leqno 2)
$$
$$
v(\lambda) + v(a) + v(x) \geq 0;
\leqno 3)
$$
$$
v(\lambda) + v(b) \geq 0;
\leqno 4)
$$
and
$$
v(\lambda) \geq 0.
\leqno 5)
$$
Eliminating $v(\lambda)$ from these, we arrive at the following
system of inequalities, together with the parity condition
$$
v(a) + v(b) + v(1-x) \, \equiv \, 0 \, (\bmod \, 2):
$$
$$
v(a) + v(x) - v(1-x) \geq 0; \leqno 6)
$$
$$
v(b) - v(1-x) \geq 0;
\leqno 7)
$$
$$
v(x) - v(1-x) \geq 0;
\leqno 8)
$$
$$
v(1-x) \leq 0;
\leqno 9)
$$
$$
v(a) + v(1-x) \leq 0;
\leqno 10)
$$
$$
v(b) + v(1-x) -v(x) \leq 0;
\leqno 11)
$$
$$
v(a) + v(b) + v(1-x) \leq 0, \leqno 12)
$$
and
$$
v(b) \geq v(a)+v(1-x). \leqno 13)
$$
To explain, $6)$ comes from $2) + 3) - 1)$, while $7)$ comes from
$2) + 4) - 1)$, $8)$ from $3) + 4) - 1)$, $9)$ from $1) - 2) - 5)$,
$10)$ from $1) - 4) - 5)$, $11)$ from $1) - 3) -5)$, $12)$ from $1)
- 2 \times 5)$, and $13)$ from $2\times 4)-1)$. The assertions of
the Lemma now follows easily.

\qed

\medskip

The first consequence is that $I_v(x)$ vanishes if $v(1-x)$ is
positive. Now let $v(1-x)$ be zero (resp. negative). Then $v(x) \geq
0$ (resp. $v(x) = v(1-x)$) and the inequalities (iii), (iv) of the
Lemma simplify to yield the conditions
$$
0 \leq v(b) \leq v(x), \quad {\rm and} \quad -v(x) \leq v(a) \leq
-v(b)
$$
(resp.
$$
v(x) \leq v(b) \leq 0 \quad {\rm and} \quad 1 \leq v(a) \leq -v(x).)
$$
The assertion of the Proposition now follows in the case $v = q \notin S$.

\bigskip

(b) \, Here $v=q$ is a prime where $\chi$ ramifies. If the conductor
of $\chi_v$ is $q^c$, with $c\geq 1$, $f_v$ is, by definition,
$g(\chi_v)^{-1}$ times a $\chi$-weighted sum of the characteristic
functions $f_{z,v}$ of $\left(
\begin{array}{cc}
1 & z \\
0 & 1
\end{array}
\right) K_vZ_v$, where $z$ has valuation $-r$ with $0\leq r \leq m$.
We need the following:

\medskip

\begin{lemma} \, $f_{z,v}\left(
\begin{array}{cc}
ab & ax \\
b & 1
\end{array}
\right)$ is non-zero only if the following conditions are
simultaneously satisfied:
\begin{enumerate}
\item[(i)] $v(1-x) \leq r$;
\item[(ii)] $v(x) = v(1-x)$ when $v(1-x) < 0$, and $v(x)=0$ when $v(1-x)>0$;
\item[(iii)] $v(a) \leq$ min$\{-v(1-x), -v(b)-v(1-x)\}$, $v(a+z)+v(ax+z)-v(a) \geq v(1-x)$;
\item[(iv)] max$\{v(1-x)+v(a), v(1-x)-v(a+z)+v(a)\} \leq v(b) \leq v(x)-v(1-x)$;
\item[(v)] $v(a) + v(b) \equiv v(1-x)$ \, $($mod $2)$.
\end{enumerate}
\end{lemma}

\medskip

{\it Proof of Lemma}. \, The statement for $r=0$ is just Lemma 5. So
we may assume that $0 < r \leq c$. Then
$$
\left(
\begin{array}{cc}
1 & z \\
0 & 1
\end{array}
\right)\left(
\begin{array}{cc}
ab & ax \\
b & 1
\end{array}
\right) \, = \, \left(
\begin{array}{cc}
ab+bz & ax+z \\
b & 1
\end{array}
\right)
$$
So for $f_{r,v}\left(
\begin{array}{cc}
ab & ax \\
b & 1
\end{array}
\right)$ to be non-zero, it is necessary and sufficient that there
exists a $\lambda$ in $\Q_v^\ast$ such that $\lambda (a+z)b, \lambda
(ax+z), \lambda b, \lambda$ are in $\Z_v$ {\it and} $\lambda^2
ab(1-x)$ is in $\Z_v^\ast$. In other words, one must have
$$
2v(\lambda) + v(a) + v(b) + v(1-x) = 0; \leqno 1')
$$
$$
v(\lambda) + v(a+z) + v(b) \geq 0; \leqno 2')
$$
$$
v(\lambda) + v(ax+z)) \geq 0; \leqno 3')
$$
$$
v(\lambda) + v(b) \geq 0; \leqno 4')
$$
and
$$
v(\lambda) \geq 0. \leqno 5')
$$
Eliminating $v(\lambda)$ from these, we arrive at the following
system of inequalities, together with the parity condition
$$
v(a) + v(b) + v(1-x) \, \equiv \, 0 \, (\bmod \, 2):
$$
$$
v(a+z)+v(ax+z)-v(a) - v(1-x) \geq 0; \leqno 6')
$$
$$
v(b) +v(a+z)-v(a)- v(1-x) \geq 0; \leqno 7')
$$
$$
v(ax+z) -v(a)- v(1-x) \geq 0; \leqno 8')
$$
$$
v(1-x) \leq v(a+z)-v(a); \leqno 9')
$$
$$
v(a) + v(1-x) \leq 0; \leqno 10')
$$
$$
v(b) + v(1-x) -v(x) \leq 0; \leqno 11')
$$
$$
v(a) + v(b) + v(1-x) \leq 0; \leqno 12')
$$
and
$$
v(b) \geq v(a)+v(1-x). \leqno 13')
$$
To explain, $6')$ comes from $2') + 3') - 1')$, while $7')$ comes
from $2') + 4') - 1')$, $8')$ from $3') + 4') - 1')$, $9')$ from
$1') - 2') - 5')$, $10')$ from $1') - 4') - 5')$, $11')$ from $1') -
3') -5')$, $12')$ from $1') - 2 \times 5')$, and $13')$ from
$2\times 4')-1')$. If $v(a)\ne v(z)$, then
$v(a+z)=$min$\{v(a),v(z)\}$, and so $v(a+z)-v(a)\leq 0$, and by
$9')$, $v(1-x)\leq 0$. On the other hand, if $v(a)=v(z)=-r$, we get
by $10')$ that $v(1-x)\leq r$. This gives part (i) of the Lemma. The
remaining assertions follow easily.

\qed

\medskip

So the first consequence is that $I_v(x)$ vanishes if $v(1-x)$ is
greater than $r$. We get from the Lemma:
$$
v(1-x)-r \leq v(b)\leq v(x)-v(1-x).
$$
This is clear when $v(a)=-r$, and if not, $v(a+z)-v(a) \leq 0$,
which even gives $v(1-x)\leq v(b)$. Next observe that when $v(1-x)$
is $0$, resp. $<0$, resp. $>0$, we have $v(x) \geq 0$, resp.
$v(x)=v(1-x)$, resp. $v(x) = 0$. We obtain
$$
-r \leq v(b) \leq v(x) \text{and} -v(x) \leq v(a) \leq 0,
$$
resp.
$$
v(x)-r \leq v(b) \leq 0 \text{and} 0 \leq v(a) \leq -v(x) \text{or}
v(a)=-r-v(x),
$$
resp.
$$
-r<v(b)<0 \text{and} v(a)=-r
$$
The assertion of the Proposition now follows in the case $v = q \in
S^\prime$, once we recall that $\vert g(\chi_v)\vert =q^{c/2}$.

\bigskip

(c) \, Let $v = N$. Here by construction, $f_N$ is the
characteristic function of $Z_N K_0(N)$ divided by $V_N$, the volume
of $Z_N K_0(N)/Z_N$. So for $f_N(\left(
\begin{array}{cc}
ab & ax \\
b & 1
\end{array}
\right) )$ to be non-zero, it is necessary and sufficient that there
exists a $\lambda$ in $\Q_N^\ast$ such that $\lambda ab, \lambda ax,
\lambda$ are in $\Z_N$, $\lambda b$ is in $N\Z_N$ {\it and}
$\lambda^2 ab(1-x)$ is in $\Z_N^\ast$. The conditions 1) through 5)
above (in the proof of part (a)) remain in force except for 4),
which gets replaced by
$$
v_N(\lambda) + v_N(b) \geq 1. \leqno 4'')
$$
The parity condition is the same as in (a). The only change in the
inequalities $6) - 12)$ is that $7)$ (resp. $8)$, resp. $10)$) gets
replaced by
$$
v_N(b) - v_N(1-x) \geq 1; \leqno 7'')
$$
$$
v(x) - v(1-x) \geq 1; \leqno 8'')
$$
and
$$
v(a) + v(1-x) \leq -1; \leqno 11'')
$$
We again get the vanishing of $I_N(x)$ when $v_N(1-x) > 0$.
Moreover, when $v_N(1-x) < 0$, we are forced to have $v_N(x) =
v_N(1-x)$, which contradicts $7')$; thus $I_N(x)$ vanishes in this
case as well. It remains only to consider when $v_N(1-x) = 0$, in
which case $v_N(x) \geq 1$. The asserted expression for $I_N(x)$
follows easily.

\bigskip

(d) \, Of course, when $r=r'=0$, $f_p$ is just the characteristic
function of $Z_pK_p$, which was treated in detail in part (a). For
general $r, r'$ with $r\geq r'$, for $f_p\left(
\begin{array}{cc}
ab & ax \\
b & 1
\end{array}
\right)$ to be non-zero, it is necessary and sufficient that there
exists a $\lambda$ in $\Q_v^\ast$ such that
$$
\left(
\begin{array}{cc}
\lambda ab & \lambda ax \\
\lambda b & \lambda
\end{array}
\right) \, \in \, K_p\left(
\begin{array}{cc}
p^r & 0 \\
0 & p^{r'}
\end{array}
\right)K_p.
$$
We get analogues of conditions $1)$ through $12)$ of part (a),
except for replacing the zeros on the right of those equations. In
$1)$ and $9)$, replace $0$ by $r+r'$, while in $2)$ through $4)$,
replace $0$ by $r$. The assertions of part d) now follow easily.

\bigskip

(e) \, Let $v = \infty$. By hypothesis, $\chi_\infty$ is $sgn$, the
sign character. By the definition of $f_\infty$, we have
$$
f_\infty\left(\begin{array}{cc}
ab & ax \\
b & 1
\end{array}
\right) \, = \, \frac{(ab(1-x))^{k/2}}{(ax-b +i(ab-1))^k} \quad
({\rm resp.} 0)
$$
if $ab(1-x)$ is $> 0$ (resp. $< 0$).

Suppose $1-x > 0$. Then the integral is over the first and the third
quadrants. Changing variables in the third quadrant and rearranging,
we get
$$
I_\infty(x) \, = \, (1-x)^{k/2} \left( I'_\infty(-1,1,1)
-(-1)^{k}I'_\infty(-1,-1,1)\right),
$$
where for $\epsilon, \delta, \nu \in \{\pm 1\}$,
$$
I'_\infty(\epsilon, \delta, \nu) \, = \,
\int_0^\infty\int_0^\infty \frac{a^{k/2 - s_1 - 1} b^{k/2 +s_2 -1}
da db} {(ax +\epsilon b +\delta i(ab +\nu))^k}.
$$
Similarly, when $1-x$ is negative, we have
$$
I_\infty(x) \, = \, (x-1)^{k/2} \left(I'_\infty(1,-1,-1)
-(-1)^{k}I'_\infty(1,1,-1)\right).
$$
Set
$$
\rho = k/2 - s_1, \quad {\rm and} \quad \sigma = k/2 + s_2.
$$
Then the assertion of part (c) will follow once we establish the
following, which was shown to us by Nathaniel Grossman.

\medskip

\begin{lemma} Suppose $k > \Re \rho > 0$ and $k > \Re \sigma > 0$. Then we have
$$
I'_\infty(\epsilon, \delta, \nu) \, = \, I_\infty(\epsilon,
\delta, \nu).
$$
This holds in the complex $x$-plane with the negative $x$-axis cut out.
\end{lemma}

\medskip

Recall the definition of $I_\infty(\epsilon, \delta, \nu)$ from
the statement of part (c) of the Proposition.

\medskip

\proof \, Fix $\epsilon, \delta, \nu$ and simply write $I'$ for
$I'_\infty(\epsilon, \delta, \nu)$. We have
\begin{eqnarray*}
I' &=& \int_0^\infty a^{\rho-1} da \int_0^\infty \frac{b^{\sigma-1}db}
{\left((ax+\delta \nu i)+b(\delta ia +\epsilon)\right)^k} \\
&=& \int_0^\infty \frac{a^{\rho-1}da}{\delta ia + \epsilon)^k} \int_0^\infty
\frac{b^{\sigma-1}db}{\left(\left(\frac{ax+\delta\nu i}{\delta ia+\epsilon}\right)
+b\right)^k }\\
&=& \int_0^\infty \frac{a^{\rho-1}\left(\frac{ax+\delta \nu i}
{\delta ia+\epsilon}\right)^\sigma da}{(\delta ia +\epsilon)^k
\left(\frac{ax+\delta \nu i}{\delta ia+\epsilon}\right)^k} \int_L \frac{b_1^{\sigma-1}db}{(1+b_1)^k},
\end{eqnarray*}
where $L$ denotes the half-line defined by the positive real
multiples of $\frac{\delta ia + \epsilon}{ax+\delta \nu i}$.

Since $\Re\lambda > \Re \sigma > 0$ by hypothesis, the ray of
integration may be rotated back to $\R_+$, and the inner
$b_1$-integral becomes
$$
\int_0^\infty \frac{b_1^{\sigma-1}db_1}{(1+b_1)^k} \, = \, B(\sigma, k-\sigma).
$$
Thus
$$
I' \, = \, \frac{B(\sigma, k-\sigma)}{x^{k-\sigma}(\delta i)^\sigma} J,
$$
where
$$
J \, = \, \int_0^\infty \frac{a^{\rho -1}da}{(a+\delta\nu i/x)^{k-\sigma}
(a-i\delta \epsilon)^\sigma}.
$$
Since by assumption $k > \Re \rho > 0$, the path of integration
defining $J$ can be rotated to the positive imaginary axis, giving
\begin{eqnarray*}
J &=& \int_0^{i\infty} \frac{a^{\rho-1}da}{(a+\delta \nu i/x)^{k-\sigma}(a-i\delta \epsilon)^\sigma} \\
&=& \int_0^\infty \frac{(ic)^{\rho-1}d(ic)}{(ic+\delta \nu i/x)^{k-\sigma}(ic-i\delta \epsilon)^\sigma} \\
&=& i^{\rho-\sigma} \int_0^\infty \frac{c^{\rho-1}dc}{(c+\delta \nu /x)^{k-\sigma}(c-\delta \epsilon)^\sigma}.
\end{eqnarray*}

\medskip

\noindent
{\it Case (i)}: \, $\delta \epsilon = -1$.

Put $u = \frac{c}{1+c}$ so that $c = \frac{u}{1-u}$ and $dc =
\frac{du}{(1-u)^2}$. We get
$$
J \, = \, i^{\rho-\sigma}(\delta \nu)^{k-\sigma} x^{k -\sigma} \int_0^1
\frac{u^{\rho-1}(1-u)^{k-\rho-1}du}{(1-u(1-\delta \nu x))^{k-\sigma}}.
$$
Now we appeal to the following well known integral representation
for the hypergeometric function ([A-S], formula (15.3.1), page 558):
$$
F(a,b;c;z) \, = \, \frac{\Gamma(c)}{\Gamma(b)\Gamma(c-b)}\int_0^\infty t^{b-1}(1-t)^{c-b-1}(1-tz)^{-a} dt,
$$
for $\Re c > \Re b > 0$.

The assertion follows by putting these together.

\medskip

\noindent
{\it Case (ii)}: \, $\delta \epsilon = 1$.

In this case,
$$
J \, = \, i^{\rho-\sigma} \int_0^\infty \frac{c^{\rho-1}dc}{(c+\delta \nu/x)^{k-\sigma}(c-1)^\sigma}.
$$
We put $c = \frac{u}{u-1}$, so that $u = \frac{c}{c-1}$ and $dc = \frac{-du}{(1-u)^2}$. We get
$$
J \, = \, i^{\rho-\sigma}(\delta \nu)^{k-\sigma}x^{k-\sigma}(-1)^{\rho-\sigma-1}
\int_0^\infty \frac{u^{\rho-1}(1-u)^{k-\rho-1}du}{(1-u(1+\delta\nu x))^{k-\sigma}}.
$$
The assertion follows by appealing once again to the integral
representation of the hypergeometric function and combining like
terms.

\qed

\vskip 0.2in

\section{\bf A bound for the sum of regular terms}

\medskip

Put
$$
I_{\rm reg}(f): = \sum_{x \in \Q^\ast-\{1\}} I(x,f).
$$

\medskip

\begin{prop} \, Let $\varepsilon >0$. Then
$$
I_{\rm reg}(f)\, \leq \, \frac{C}{N^{k/2-\varepsilon}},
$$
for a positive constant $C$.
\end{prop}

\medskip

\begin{proof}. \, Put $t=\frac{1}{1-x}$. Since $x \to t$ is an automorphism of $\Q^\ast-\{1\}$,
$$
I_{\rm reg}(f): = \sum_{t \in \Q^\ast-\{1\}} I(x,f),
$$
with $x=\frac{t-1}{t}$. By Proposition 2.4, there is a positive
integer $M=Dp^r$ such that for any finite place $v$,
$$
v(t)<-v(M) \, \implies \, I_v(x,f) = 0.
$$
In fact that Proposition says that for $v=N$, $I_v(x,f) \ne 0
\implies N \mid (t-1)$. Also, $I_\infty(x,f)=0$ unless $x>0$.
Putting these together, we see that $I(x,f)$ is zero unless $Mt$
is an integer $\ne 0,1$, which is divisible by $N$. Thus
$$
I_{\rm reg}(f): = \sum_{n\ne 0,1, \, {N\vert (n-M)}}
I(\frac{n-M}{n},f).
$$

\begin{lemma} \, For any $n \ne 0$,
$$
I(\frac{n-M}{n},f) \, << \, \frac{1}{n^{k/2-\varepsilon}},
$$
with the implied constant independent of $n$, \end{lemma}

\medskip

\begin{proof} \, By Prop. 2.4, for any (finite) prime $q$,
$I_q(\frac{n-M}{n},f)$ is $1$ for $q \nmid n(n-M)$, and if $q\mid
n(n-M)$,
$$
I_q(\frac{n-M}{n},f) \, \leq M v_q((n-M)/n)^2.
$$
We {\it claim} that the function
$$
g(n): = \, \prod_q v_q(n)
$$
satisfies, for every $\varepsilon>0$, the bound
$$
g(n) \, \leq \, Cn^\varepsilon,
$$
for a constant $C>0$ independent of $n$.

To see this, first note that $g(n)$ is multiplicative in $n$
(though not strictly). Fix any $\varepsilon>0$. Then $\exists A>0$
such that for any prime $q$ and positive integer $a>0$,
$$
a \leq p^{a\varepsilon} \text{if either}  a > M  \text{or}  p > M.
$$
Given $n >0$ with unique factorization $n= \prod_{j\in J}
q_j^{a_j}$, where the $q_j$ are primes and $a_j > 0$, we may write
$n=mk$, where $m$ has no prime divisors  $q_j  <  A$ with
exponents $a_j < A$. It follows immediately that $g(m) \leq
nm^\varepsilon$,  and $k$  is a product over primes $\leq A$ with
exponents at most $A$.  We obtain
$$
g(n) = g(m)g(k) \leq  Cn^e
$$
where  $C$  is the maximum of  $g(k)$  for the finitely many
choices of $k$.

Hence the {\it claim}.

\medskip

It follows that for any $\varepsilon
> 0$,
$$
\prod\limits_{q\mid n(n-M)}v_q((n-M)/n)^2 \, << \, n^\varepsilon.
$$
Moreover, thanks to Prop. 2.4, we have
$$
I_\infty(\frac{n-M}{n},f) \, << \, \frac{1}{n^{k/2}}.
$$
In both estimates, the implied constants depends on $M$, but {\it
not} on $N$. This proves the Lemma.

\end{proof}

\medskip

Consequently,
$$
I_{\rm reg}(f) \, << \, \sum_{n\ne 0,1, \, {N\vert n-M}}
\frac{1}{n^{k/2-\varepsilon}}.
$$
But for any real $u >1$,
$$
\sum_{n\ne 0,1, \, {N\vert n-M}} \frac{1}{n^u} \, = \, \sum_{m\geq
1}\frac{1}{(mN+M)^u} \, << \, \frac{1}{N^u}.
$$
The Proposition now follows.
\end{proof}

\vskip 0.2in

\section{\bf The spectral side}

\bigskip

Let $I(f)$ be the limit, as $s_1, s_2 \to 0$, of the integral of the
kernel of the test function over the square of the torus $\tilde T$
against the quasi-character
$(\chi\vert\cdot\vert^{s_1},\vert\cdot\vert^{s_2})$ with $\chi$ odd.
Let $\rho$ denote the right regular representation of $\tilde G(\A)$
on the space of cusp forms with trivial central character.

\medskip

Let $\{\varphi_j\}$ denote an orthogonal basis of cusp forms of
level $N$ and weight $k\geq 2$. They can be chosen to be Hecke
eigenforms, i.e., generate an admissible representation
$\pi_j=\pi_{j,\infty}\otimes\pi_{j,0}$ of GL$(2, \A)$ satisfying
$\pi_{j,\infty} \simeq {\mathcal D}_k$ and $\pi_{j,0}^{K_0(N)} \ne
0$. For the subspace of newforms of level $N$, we may take
$\varphi_j$ to be newforms. If $k<12$, then we have only newforms.

\medskip

Write $f=f_\infty\otimes f_0$, where $f_0$ denotes the finite part
of $f$. By construction $f_\infty$ is $d({\mathcal D}_k)$ times the
complex conjugate of the matrix coefficient $h_\infty$, say, of
$\mathcal D$. When we replace the Haar measure $dg$  by  $c dg$, the
formal degree gets multiplied by $c^{-1}$. The matrix coefficient $h
_\infty$  is independent of the choice of Haar measure, but the
operator $\rho_\infty(h_\infty)$ does depend  on the Haar measure.
By contrast, the function $f_\infty$  depends on the Haar measure,
but the operator  $\rho_\infty(f_\infty )$  is independent of the
Haar measure. In fact, $\rho_\infty(f_\infty )$ is a projection
operator onto the subspace of cusp forms of weight $k$. Note also
that the kernel of a fixed function  $h(g)$ does not depend on $dg$
since it's just the sum  $\sum h(g^{-1} \gamma h)$. However,  the
kernel for $\rho(f)$ depends inversely on Haar measure because the
formal degree does.

Putting these remarks together, we obtain the spectral expression
$$
K_f(x,y) \, = \, \sum_j
\frac{\rho(f_0)\varphi_j(x)\overline{\varphi_j(y)}}{\langle
\varphi_j,\varphi_j\rangle},
$$
where the sum is over an orthogonal basis of cusp forms of weight
$k$ and level $N$. On each $\varphi_j$, $\rho(f)$ acts by
$\pi_j(f)$. There is an adelic scalar product $\langle \cdot,
\cdot\rangle_\A$ on the space of cusp forms, given by
$$
\langle \varphi, \psi\rangle_\A \, = \,
\int\limits_{Z(\A)G(\Q)\backslash G(\A)}
\varphi(g)\overline{\psi(g)} dg,
$$
where $dg$ denotes the quotient measure on $Z(\A)G(\Q)\backslash
G(\A)$ induced by the Haar measure on $G(\A)$, normalized so that
the quotient space gets volume $1$. One has the identification
$$
Z(\A)G(\Q)\backslash G(\A)/K_\infty K_0(N) \, = \,
\Gamma_0(N)\backslash {\mathcal H},
$$
and if $\varphi$ is a cusp form of level $N$,
$\vert\varphi\vert^2$ is right invariant by $K_\infty K_0(N)$.
Consequently, $\langle \varphi, \varphi\rangle_\A$ is $V_N$ times
the Petersson scalar product $\langle \varphi, \psi\rangle$ on the
classical modular forms of level $N$. By abuse of notation, we are
using the same symbol for the adelic cusp form and the classical
modular form it defines.

\medskip

\begin{prop} \, We have
$$
I(f) \, = \, \frac{1}{4\pi V_N}\sum_j \frac{L(1/2,
\varphi_j)L(1/2,\varphi_j\otimes \chi)}{\langle
\varphi_j,\varphi_j\rangle}\hat f_p(\varphi_j),
$$
where $\hat f_p(\varphi_j)$ denotes a local expression depending
only on $f_p$ and $\varphi_j$. For a suitable local function
$f_{p,0}$,
$$
\hat f_{p,0}(\varphi_j) \, = \, a_p(\varphi_j).
$$
\end{prop}

\medskip

\begin{proof}. \, Using the expression above for the kernel in
conjunction with the definition of $I(f)$, we obtain
$$
I(f) \, = \, \lim_{s_1\to 0, s_2\to 0} \, \sum_j P_j(f,\chi,
s_1)Q_j(f, s_2),
$$
where $P_j, Q_j$ are period integrals given by
$$
P_j(f,\chi, s_1) \, = \, \int_{\Q^\ast\backslash\A^\ast}
\rho(f)\varphi_j\left(
\begin{array}{cc}
a & 0 \\
0 & 1
\end{array}
\right) \chi(a)\vert a\vert^{s_1} d^\ast a
$$
and
$$
Q_j(f,s_2) \, = \, \int_{\Q^\ast\backslash\A^\ast}
\overline\varphi_j\left(
\begin{array}{cc}
b & 0 \\
0 & 1
\end{array}
\right) \vert a\vert^{s_2} d^\ast b .
$$
If $W_j$ denotes the Whittaker function of $\varphi_j$, there is a
well known Fourier expansion (for $g \in {\rm GL}(2,\A)$)
$$
\varphi_j(g) \, = \, \sum_{t \in \Q^\ast} \, W_j\left[\left(
\begin{array}{cc}
t & 0 \\
0 & 1
\end{array}
\right)g\right].
$$
Since $\chi$ and the adelic absolute value $\vert \cdot \vert$ are
$1$ on $\Q^\ast$, we can unfold the expression for $P_j$ as
$$
P_j(f,\chi, s_1) \, = \, \int_{\A^\ast} \rho(f)W_j\left(
\begin{array}{cc}
a & 0 \\
0 & 1
\end{array}
\right) \chi(a)\vert a\vert^{s_1} d^\ast a
$$
Then from the factorizability of $W_j$, $f$ and $\chi$, we get the
product expansion
$$
P_j(f,\chi, s_1) \, = \, \prod_v P_{j,v}(f,\chi,s_1),
$$
where $v$ runs over all the places of $\Q$, with
$$
P_{j,v}(f,\chi,s_1) \, = \, \int_{\Q_v^\ast} \rho(f_v)W_{j,v}\left(
\begin{array}{cc}
a & 0 \\
0 & 1
\end{array}
\right) \chi_v(a)\vert a\vert^{s_1} d^\ast a .
$$

\medskip

For $v \ne p$, $f_v$ is chosen such that $\rho(f_v)W_{j,v}$ times
$\chi_v\circ {\rm det}$ is the new vector of $\pi_{j,v}\otimes
\chi_v$, so that by Jacquet-Langlands,
$$
P_{j,v}(f,\chi,s_1) \, = \, L(s_1+1/2,\pi_{j,v}\otimes \chi_v).
$$
This is clear at any $v$ where $\chi$ is unramified. Let us indicate
the reason at a $v=q$ where $\chi$ is ramified, say of conductor
$q^m$. It follows from the definition of $f_v$, and the
transformation property of the Whittaker function under the left
translation by the maximal unipotent subgroup, that
$$
P_{j,v}(f,\chi,s_1) \, = \, \frac{1}{g(\chi_v)}\sum_{n \in \Z}
W_v\begin{pmatrix}\varpi^n & 0\\0 & 1\end{pmatrix}
g(\chi_v,\psi_{\varpi^{n-m}})q^{-n(s_1+1/2)},
$$
where $\psi_{t}(x)=\psi(tx)$ and $\psi$ the additive character
defined as the composite
$$
\Q_q \rightarrow \Q_q/\Z_q \rightarrow \Q/\Z \rightarrow S^1,
$$
with the last arrow on the right being $x \to e^{2\pi i x}$. Note
that $g(\chi_v)$ is just $g(\chi_v,\psi_{\varpi^{-m}})$. Since $W_v$
has been chosen to be the spherical vector giving the right
$L$-factor of $\pi_v$, we have
$$
W_v\begin{pmatrix}\varpi^n & 0\\0 & 1\end{pmatrix} \, = \,
\delta_{n\geq 0} \, \,
q^{n/2}\left(\frac{\alpha^{n+1}-\beta^{n+1}}{\alpha-\beta}\right),
$$
where $\delta_{n\geq 0}$ is $1$ if $n \geq 0$ and is zero otherwise,
and ${\it diag}\{\alpha,\beta\}$ is the conjugacy class attached to
$\pi_v$. (In other words, $L(s, \pi_v)$ is the inverse of $(1-\alpha
q^{-s})(1-\beta q^{-s})$.) So we may take $n$ to be non-negative
from here on. On the other hand, it is well known that
$g(\chi_v,\psi_{\varpi^r})=0$ is $r<m$. Taking $r=m-n$ we get
vanishing for $n >0$. Thus $n=0$, and we get (for any $s_1$)
$$
P_{j,v}(f,\chi,s_1) \, = \, 1,
$$
To finish, note that
$$
L(s,\pi_v\otimes\chi_v) \, = \, 1,
$$
because $\pi_v$ is unramified, while $\chi_v$ is ramified.

\bigskip

Similarly, $\forall v \ne p$,
$$
Q_j(f,s_2) \, = \, \prod_v Q_{j,v}(f,s_2),
$$
with
$$
Q_{j,v}(f,s_2) \, = \, L(\overline s_2+1/2,\pi_{j,v}), \\
\forall \, v .
$$
The main assertion of the Proposition now follows, after letting
$s_1$ and $s_2$ go to $0$, which is justified because the
expressions are known to have analytic continuations, with
$$
\hat f_p(\varphi_j) \, = \,
\frac{P_{j,v}(f,\chi,0)}{L(1/2,\pi_{j,p}\otimes\chi_p)}.
$$
Finally, one knows that there is a Hecke function $f_{p,0}$ such
that
$$
\pi_j(f_{p,0})\varphi_j \, = \, a_p(\pi_j)\varphi_j,
$$
and this gives the final assertion.

\end{proof}

\vskip 0.2in

\section{\bf Proof of the Main Theorem sans measure}

\bigskip

Now we will exploit the equality of the geometric and spectral
sides of the relative trace formula at hand, for forms of even
weight $k > 2$ and prime level $N$. Continue to assume that $\chi$
is a quadratic character of discriminant $D <0$.

The sum of the regular terms on the geometric side was estimated in
section 3. There are six singular terms, and four of them vanish. As
$s_j\to 0$, for $j=1,2$, each of the remaining pair of singular
terms, which are the dominant ones, goes to $L(1,\chi)$. In view of
the expression for the spectral side (cf. section 4), we then obtain
the following (for any $\varepsilon
> 0$):
$$
\frac{1}{4\pi}\sum_{\varphi \in \mathcal F_{N,k}}
\frac{1}{(\varphi,\varphi)} L(\frac{1}{2}, \varphi \otimes
\chi)L(\frac{1}{2}, \varphi) \, = \, 2c_kL(1, \chi) +
O(N^{-k/2+\varepsilon}),
$$
where $\mathcal F_{N,k}$ is an orthogonal basis of holomorphic,
Hecke eigencusp forms of weight $k$ and level $N$, and $c_k$ is as
in the statement of Theorem A.

\medskip

Now let $\mathcal F_N(k)^{\rm new}$ be the subspace of newforms of
level $N$ and even weight $k$. The sum on the left of the formula
above runs over newforms of level $N$ and oldforms which
necessarily come from level $1$. Since $\chi(-1)=-1$, we see that
the sign of the functional equation of $L(s,g\otimes \chi)L(s,g)$
is $-1$ for any form of level $1$, implying that $L(1/2,g \otimes
\chi)L(1/2,g)$ is zero for any such $g$. There are two oldforms of
level $N$ associated to $g$, which, in classical language, are
$g(z)$ and $g_N(z):=g(Nz)$. It is immediate that their Mellin
transforms are related by $L(s,g_N) = N^{1-s}L(s,g)$. It follows
that $L(1/2,g_N\otimes \chi)L(1/2,g_N)$ is also $0$.

Now we have proved Theorem A without the assertion about the
measure.

\vskip 0.2in

\section{\bf The measure}

\medskip

Fix $v=p$ not dividing $ND$, so that $\pi$ and $\chi$ are both
unramified at $p$.

Let $f_n$ be the characteristic function of \ds{ZK\left(
\begin{array}{cc}
  p^n & 0 \\
  0 & 1 \\
\end{array}
\right)K} for $n= 0, 1, 2,...$ and let $\varphi_n(s)$  be the Satake
transform of $f_n$. Then $\varphi_n(s)$ is a symmetric Laurent
polynomial in $p^s$ and $p^{-s}$ with the property:
$$
\textrm{Trace}(\pi_s(f_n)) = \varphi_n(s)
$$
where $\pi_s$ is the principal series representation unitarily
induced from the character
$$
\left(\begin{array}{cc}
 \alpha & 0 \\
  0 & \beta \\
\end{array}
\right) \to \left|\frac{\alpha}{\beta}\right|^s
$$

\begin{lemma}
$\varphi_0(s)=1$ and for $n\ge 1$,
$$
\varphi_n(s) =
p^{n/2}\Big[p^{ns}+p^{-ns}+\left(1-\frac1p\right)\Big(p^{(n-2)s}
+p^{(n-4)s}+\cdots +p^{-(n-4)s})+p^{-(n-2)s}\Big) \Big]
$$
\end{lemma}

\begin{proof}
The coset  \ds{K\left(
\begin{array}{cc}
  p^n & 0 \\
  0 & 1 \\
\end{array}
\right)K} is equal to a union of single cosets  $gK$ where $g$
ranges over the following elements:
\begin{alignat*}{2}
&\left(
\begin{array}{cc}
  p^n & t \\
  0 & 1 \\
\end{array}
\right)&\qquad\qquad &t\in \mathcal O/(p^n)\\
&\left(
\begin{array}{cc}
  p^{n-k} & t \\
  0 & p^k \\
\end{array}
\right)&\qquad\qquad &t\in \Big(\mathcal O/p^{n-k}\Big)^*\quad k = 1, ..., n\\
\end{alignat*}
Using this we easily calculate the action of $f_n$ on the
principal series $\pi_s$ and we find the above formula.
\end{proof}

We want to find the Satake transform of the following linear
functional for bi-$K$-invariant $f$:
\begin{eqnarray*}
T(f) = I_{v}(n^{+})  \ &=&\iint_{\Q_{v}^{*}\times \Q_{v}^{*}}\
f(\left(
\begin{array}{cc}
b & a \\
0 & 1
\end{array}
\right) )\ \chi _{v}|a|^{-s_{1}-s_{2}}\chi _{2,v}(b)|b|^{s_{2}}\
d^{*}a\ d^{*}b.
\end{eqnarray*}
In fact, we will take $\chi_2$ trivial, $s_2 = 0$, and set
$s=s_1$. We assume $\chi$ is unramified, so we have:
\begin{eqnarray*}
T(f) = I_{v}(n^{+})  \ &=&\iint_{\Q_{v}^{*}\times \Q_{v}^{*}}\
f(\left(
\begin{array}{cc}
b & a \\
0 & 1
\end{array}
\right) )\ \chi(a)|a|^{-s} \ d^{*}a\ d^{*}b.
\end{eqnarray*}

We first evaluate  $T(f_n)$ for $n\ge 0$. We first note that if
$$
\left(
\begin{array}{cc}
  b & a \\
  0 & 1 \\
\end{array}
\right)\,\,\in\,\,ZK\left(
\begin{array}{cc}
  p^n & 0 \\
  0 & 1 \\
\end{array}
\right)K
$$
then
$$
\left(
\begin{array}{cc}
  b & a \\
  0 & 1 \\
\end{array}
\right) \in K\left(
\begin{array}{cc}
  p^{n+\alpha} & 0 \\
  0 &  p^{\alpha} \\
\end{array}
\right)K
$$
for some $\alpha$. This is not possible unless  $\alpha\le 0$ and
$b=p^{n+2\alpha}u$ for some unit $u$. Furthermore,  the gcd of
$\{a, p^{n+2\alpha}, 1\}$ is $p^{\alpha}$, so we must have
$n+2\alpha\ge \alpha$, that is, $\alpha\ge -n$.

If  $0>\alpha> -n$, then the gcd condition forces $a\in
p^{\alpha}\Z_p^*$. If $\alpha=0$, then we just have $a\in\Z_p$ and
if $\alpha=-n$, we have $\alpha\in p^{-n}\Z_p$.

Set $\delta = \chi(p)$. For $n=0$, we have:
\begin{align*}
 T(f_0) =  \int_{u\in \Z_{p}^{*}}\int_{a\in\Z_{p}}  \chi(a)|a|^{-s} \ d^{*}u\
 d^{*}a
 = L_p(-s,\chi)
\end{align*}
For $n>0$, $T(f_n)$ is a sum of three terms:
$$
I = \int_{u\in \Z_{p}^{*}}\int_{a\in\Z_{p}}  f(\left(
\begin{array}{cc}
p^nu & a \\
0 & 1
\end{array}
\right) )\chi(a)|a|^{-s} \ d^{*}u\ d^{*}a
 = L_p(-s,\chi)
$$
and
$$
II = \int_{u\in \Z_{p}^{*}}\int_{a\in p^{-n}\Z_{p}}  f(\left(
\begin{array}{cc}
p^{-n}u & a \\
0 & 1
\end{array}
\right) )\chi(a)|a|^{-s} \ d^{*}u\ d^{*}a
 = \delta^{-n}p^{-ns}L_p(-s,\chi)
$$
\begin{align*} III &= \sum_{\alpha=1}^{n-1}\,\int_{u\in
\Z_{p}^{*}}\int_{v\in\Z_{p}^*} f(\left(
\begin{array}{cc}
p^{n-2\alpha}u & p^{-\alpha}v \\
0 & 1
\end{array}
\right) )\chi(p^{-\alpha})|p^{-\alpha}|^{-s} \ d^{*}u\ d^{*}v\\
&=  \sum_{\alpha=1}^{n-1}\chi(p^{-\alpha})p^{-\alpha s} =
\frac{\delta^{-1}p^{-s}-\delta^{-{n}}p^{-{n}s}}{1-\delta^{-1}p^{-s}}
=
\frac{\delta^{1-n}p^{(1-n)s}-1}{1-\delta p^{s}}\\
& = L_p(-s,\chi)(\delta^{1-n}p^{(1-n)s}-1)
\end{align*}
So we get the sum:
$$
T(f_n) = L_p(-s,\chi)(\delta^{-n}p^{-ns}+\delta^{1-n}p^{(1-n)s}) =
\delta^{-n}p^{-ns}L_p(-s,\chi)(1+\delta p^{s})
$$

Note that if $\delta = -1$, then $T(f_n)$ this approaches $0$ as
$s\to 0$ and in this, the limit is the {\it Plancherel measure} at
$p$, that is, the measure $\mu_-$ satisfying $\mu(f_0)=1$ and
$\mu(f_n)=0$ for $n>0$.

\bigskip

Now we calculate
\begin{eqnarray*}
S(f) = I_{v}(n^{-})  \ &=&\iint_{\Q_{v}^{*}\times \Q_{v}^{*}}\
f(\left(
\begin{array}{cc}
a & 0 \\
b & 1
\end{array}
\right) )\ \chi(a^{-1}b)|a|^{-s_1}\,|b|^{s_{1}+s_{2}} d^{*}a\
d^{*}b.
\end{eqnarray*}

We  evaluate  $S(f_n)$ for $n\ge 0$. We first note that if
$$
\left(
\begin{array}{cc}
  a & 0 \\
  b & 1 \\
\end{array}
\right)\,\,\in\,\,ZK\left(
\begin{array}{cc}
  p^n & 0 \\
  0 & 1 \\
\end{array}
\right)K
$$
then
$$ \left(
\begin{array}{cc}
  a & 0 \\
  b & 1 \\
\end{array}
\right) = \left(
\begin{array}{cc}
  p^nu & 0 \\
  b & 1 \\
\end{array}
\right),\,\, \left(
\begin{array}{cc}
  p^{n-2\alpha}u & 0 \\
   p^{-\alpha}v & 1 \\
\end{array}
\right)\,\, \textrm{or}\,\,\,\left(
\begin{array}{cc}
  p^{-n}u & 0 \\
   p^{-n}b & 1 \\
\end{array}
\right)
$$
where $b\in \Z_p$ and $u,v$ are units in $\Z_p^*$, and  $0\le
\alpha\le n$.

For $n=0$, we have:
\begin{align*}
 T(f_0) =  \int_{u\in \Z_{p}^{*}}\int_{b\in\Z_{p}}  \chi(b)|b|^{s_1+s_2} \ d^{*}u\
 d^{*}a
 = L_p(s_1+s_2,\chi)
\end{align*}
For $n>0$, $T(f_n)$ is a sum of three terms:
$$
I = \int_{u\in\Z_p^{*}} \int_{b\in\Z_p}
\chi(p^{-n}b)|p^n|^{-s_1}\,|b|^{s_{1}+s_{2}} d^{*}u d^{*}b =
\delta^{-n}p^{ns_1}L(s_1+s_2,\chi)
$$

\begin{align*} II &=
\sum_{\alpha=1}^{n-1}\,\,\int_{u\in\Z_p^{*}}\int_{v\in\Z_p^*}
\chi(p^{\alpha-n})|p^{n-2\alpha}|^{-s_1}\,|p^{-\alpha}|^{s_{1}+s_{2}}
d^{*}u\ d^{*}v\\
&=p^{ns_1}\delta^{-n}\sum_{\alpha=1}^{n-1}\delta^\alpha\,p^{\alpha(s_2-s_1)}
 =p^{ns_1}\delta^{-n}\frac{\delta p^{s_2-s_1} - \delta^2
p^{n(s_2-s_1)}}{1 - \delta p^{s_2-s_1}}
 =p^{ns_1}\delta^{-n}\frac{\delta p^{s_2-s_1} - \delta^2 p^{n(s_2-s_1)}}{1 - \delta p^{s_2-s_1}}\\
\end{align*}

\begin{align*} III &= \int_{u\in\Z_p^{*}}\int_{b\in\Z_p}
  \chi(b)|p^{-n}|^{-s_1}\,|p^{-n}b|^{s_{1}+s_{2}} d^{*}a\
d^{*}b =  p^{ns_2}L(s_1+s_2,\chi)
\end{align*}

When we go from $I_v(n^+)$ to $I_v(n^-)$, the main change is that
$(s_1,s_2)$ gets sent to $(-s_2,-s_1)$. We get $I_v(n^+)=I_v(n^-)$
(as $s_1, s_2 \to 0$), and we again get the Plancherel measure at
$p$ when $\chi(p)=-1$.

\medskip

We are left to analyze the case $\chi(p)=1$. We will simply write
$\mu$ for $\mu_+$. By the calculation above, we have
$$
T(f_0) = 1, \text{and} T(f_n) = \delta^n +\delta^{-n} \, \,
\forall n
>0,
$$
with $\delta=1$. It will be convenient (and more transparent) for us
to carry out the calculations for arbitrary non-zero scalar $\delta$
and at the end set $\delta=1$. We seek a function $F(s)$ of $p^s$
such that
$$
\int \varphi_n(s)F(s)\,ds = \left\{\begin{array}{cc}
  1 & n=0 \\
  \delta^n+\delta^{-n} & n\ge 1 \\
\end{array}\right.
$$
The integral is taken over $s\in i\mathbf R/(2\pi i \ln p \,\mathbf
Z)$.

\medskip

\begin{prop}
$$
F_{ev}(s) = 1 + \sum_{n=1}^\infty\, A_n p^{2ns-n}
$$
where

\begin{align*} A_n &= \delta^{2n}+\delta^{-2n} -
(p-1)\Big(\delta^{2n-2}+\delta^{2n-4}+\delta^{2n-6}+\cdots +
\delta^{-2n+4}+\delta^{-2n+2}\Big)\\ &= \delta^{2n}+\delta^{-2n} -
(p-1)\Big(\frac{\delta^{2n-1}-\delta^{-2n+1}}{\delta-\delta^{-1}}\Big)
\end{align*}

\end{prop}

\begin{proof}
 Observe that
$$
\int p^{ns}\,p^{ms}\,ds =  \left\{\begin{array}{cc}
  1 & n+m=0 \\
  0 & \textrm{otherwise} \\
\end{array}\right.
$$
Therefore \ds{\int \varphi_{2n}(s)F_{ev}(s)\,ds} is equal to
$$
A_n+\Big(1-\frac1p\Big)\Big[pA_{n-1}+p^2A_{n-2}+\cdots p^nA_0\Big]
= \delta^{2n}+\delta^{-2n}+\sum_{j=-n+1}^{n-1}c_j\delta^{2j}
$$
We must show that the coefficients $c_j$ are zero. We have:

\begin{align*}
c_j &=
-(p-1)-\Big(1-\frac1p\Big)\Big[p(p-1)+p^2(p-1)+\cdots+p^{j-1}(p-1)-p^j\Big]\\
&=-(p-1)-\Big(1-\frac1p\Big)(-p) = 0
\end{align*}
\end{proof}

\begin{prop}
$$
F_{od}(s) =  B_0p^{s-\frac12} +
p^{3s-\frac32}B_1+p^{5s-\frac52}B_2+\cdots = \sum_{n\ge 1}
B_np^{(2n+1)s-(n+\frac12)}
$$
where
\begin{align*} B_n &= \delta^{2n+1}+\delta^{-2n-1} -
(p-1)\Big(\delta^{2n-1}+\delta^{2n-3}+\delta^{2n-5}+\cdots +
\delta^{-2n+3}+\delta^{-2n+1}\Big)\\ &=
\delta^{2n+1}+\delta^{-2n-1} -
(p-1)\Big(\frac{\delta^{2n}-\delta^{-2n}}{\delta-\delta^{-1}}\Big)
\end{align*}

\end{prop}

\begin{proof}
A similar computation gives that \ds{\int
\varphi_{2n+1}(s)F_{ev}(s)\,ds} is equal to
$$
B_n+\Big(1-\frac1p\Big)\Big[pB_{n-1}+p^2B_{n-2}+\cdots p^nB_0\Big]
= \delta^{2n+1}+\delta^{-2n-1}+\sum_{j=-n+1}^{n-1}c_j\delta^{2j+1}
$$

\qed
\end{proof}

In short, we get
$$
F(s) = 1 + \sum_{n=1}^\infty\,C_np^{n(s-\frac12)}
$$
where
$$
C_n = \delta^n+\delta^{-n}-
(p-1)\Big(\frac{\delta^{n-1}-\delta^{-n+1}}{\delta-\delta^{-1}}\Big)
$$
To write $F(s)$ as a rational function, we sum the geometric
series. Let $T =p^{(s-\frac12)}$.

\begin{align*}
F(s)&= 1 + \sum_{n\ge 1}\Big[ (\delta T)^n + (\delta^{-1} T)^n -
\frac{(p-1)}{\delta-\delta^{-1}}  \Big(  \delta^{-1}(\delta T)^n +
\delta(\delta^{-1} T)^n
   \Big)
   \Big]\\
   &= 1 + \sum_{n\ge 1}\frac{(\delta
   T)^n}{\delta-\delta^{-1}}\Big[\delta-\delta^{-1}-p\delta^{-1}+\delta^{-1}\Big]
+ \sum_{n\ge 1}\frac{\delta^{-1}
   T)^n}{\delta-\delta^{-1}}\Big[\delta-\delta^{-1}+p\delta - \delta \Big]\\
   &= 1 +  \frac{\delta -p\delta^{-1}}{\delta-\delta^{-1}}
\sum_{n\ge 1}(\delta T)^n - \frac{\delta^{-1}
-p\delta}{\delta-\delta^{-1}}
\sum_{n\ge 1}(\delta^{-1} T)^n \\
  &= 1 +  \frac{\delta -p\delta^{-1}}{\delta-\delta^{-1}}
\frac{ \delta T}{1-\delta T} - \frac{\delta^{-1}
-p\delta}{\delta-\delta^{-1}}
\frac{ \delta^{-1} T}{1-\delta^{-1} T} \\
\end{align*}
Using Mathematica, we find the further simplification:
$$\boxed{
F(s) = \frac{1-pT^2}{(T-\delta)(T-\delta^{-1})} =
\frac{1-p^{2s}}{(1- \delta^{-1}p^{s-\frac12})(1-\delta
p^{s-\frac12})} }
$$

We can also replace $F(s)$ by $\Re(F(s)) = \frac12(F(s)+F(-s))$
(for $s$ purely imaginary). We find that
$$
\mu(s) = \frac12(F(s)+F(-s)) =
\frac12\frac{\Big(1-\frac1p\Big)(2-p^{2s}-p^{-2s})}{\Big(1+\frac1p-p^{-\frac12+s}-p^{-\frac12-s}\Big)^2}
$$
The linear functional is
$$
T(f) = \frac1T\int_{0}^{T} \widehat{f}(is)\,\mu(is)ds
$$
where $T = \pi /\log(p)$. Here is a graph of $\mu(s)$ for the case
$p=2$ and $\delta$ trivial.

\centerline{\includegraphics{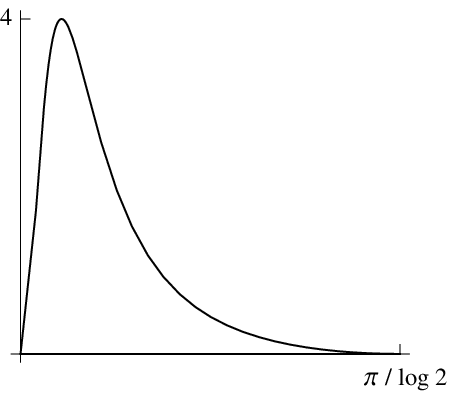}}

Now let $x = p^s+p^{-s} = 2\cos(s\log p)$. Then
$x^2=2+p^{2s}+p^{-2s}$ and we can also write
\begin{align*} \mu(s) &= \frac12(F(s)+F(-s)) =
\frac12\frac{\Big(1-\frac1p\Big)(4-x^2)}{\Big(1+\frac1p-p^{-\frac12}
x\Big)^2}\\
dx & = -2\log(p)\sin(s\log p)\,ds =
-2\log(p)\sqrt{1-(x/4)^2}\,ds=- \log(p)\sqrt{4-x^2}\,ds
\end{align*}
and thus

\begin{align*}
\frac{\log(p)}{\pi}\mu(s)ds  &=  - \Big(\frac{\log(p)}{\pi} \Big)
\frac1{2\log(p)}\frac{\Big(1-\frac1p\Big)\sqrt{4-x^2}}{\Big(1+\frac1p-p^{-\frac12}
x\Big)^2}\, dx\\
&=  -
\frac{\pi}{2}\frac{\Big(1-\frac1p\Big)\sqrt{4-x^2}}{\Big(1+\frac1p-p^{-\frac12}
x\Big)^2}\, dx\\
\end{align*}
Notice that the minus sign disappears when we write the integral
as an integral over $[-2,2]$.

\noindent Dinakar Ramakrishnan\\
Professor of Mathematics, 253-37 Caltech, Pasadena, CA 91125, USA\\
E-mail: dinakar@its.caltech.edu

\noindent Jonathan Rogawski\\
 Professor of Mathematics Department of Mathematics\\
 University of California, Los Angeles, California 90095\\
E-mail: jonr@math.ucla.edu

\end{document}